\documentclass[10pt]{amsart}
\usepackage{enumerate,amsmath,amssymb,latexsym,
amsfonts, amsthm, amscd}

\theoremstyle{plain}

\newtheorem{theorem}{Theorem}[section]

\numberwithin{equation}{section}

\newcommand{\End}{{\rm End}}
\newcommand{\IM}{{\rm Im}}
\newcommand{\ra}{\rightarrow}

\addtocounter{section}{-1}

\begin{document}

\title {The bosonic Fock representation and a generalized Shale Theorem}

\date{}

\author[P.L. Robinson]{P.L. Robinson}

\address{Department of Mathematics \\ University of Florida \\ Gainesville FL 32611  USA }

\email[]{paulr@ufl.edu}

\subjclass{} \keywords{}

\begin{abstract}

We detail a new approach to the bosonic Fock representation of a complex Hilbert space $V$: our account places the bosonic Fock space $S[V]$ between the symmetric algebra $SV$ and its full antidual $SV'$; in addition to providing a context in which arbitrary (not necessarily restricted) real symplectic automorphisms of $V$ are implemented, it offers simplified proofs of many standard results of the theory.

\end{abstract}

\maketitle


\medbreak

\section{Introduction}

Traditionally, the bosonic Fock representation of a complex Hilbert space $V$ is founded in symmetric Fock space $S[V]$: the Hilbert space completion of the symmetric algebra $SV$ relative to a canonical product. Again traditionally, the various operators of interest (such as the number operator, field operators, creators and annihilators) are initially defined on the symmetric algebra and then extended to their maximal domains in Fock space. An unfortunate aspect of this traditional approach is that these extended operators are defined implicitly rather than by explicit formulae, a circumstance that often entails the use of awkward and indirect arguments. 

\medbreak 
A celebrated theorem of Shale asserts that a symplectic automorphism $g$ of $V$ is unitarily implemented in the Fock representation on $S[V]$ if and only if the commutator $[g,i] = gi - ig$ is a Hilbert-Schmidt operator. A standard proof of this theorem involves first developing an essentially figurative expression for the corresponding displaced vacuum and then showing that the Hilbert-Schmidt condition is necessary and sufficient for this figurative expression to define an element of $S[V]$. It is reasonable to ask for a context in which such figurative expressions are strictly legitimate: a setting that accommodates displaced vacua for all symplectic automorphisms. 

\medbreak 
Our purpose in these notes is to present a new approach to the bosonic Fock representation that addresses each of the issues just mentioned. In spirit, ours is a variant of the rigged Hilbert space approach and places $S[V]$ between a suitable subspace and its antidual. In fact, we follow the simplest route: the canonical inner product embeds $SV$ in its full (purely algebraic) antidual $SV'$ comprising all antilinear functionals $SV \ra {\mathbb C}$; Fock space $S[V]$ is realized as the subspace of bounded antilinear functionals, whence the triple $SV \subset S[V] \subset SV'$. An important feature of this approach is that the antidual $SV'$ is itself a commutative associative algebra: indeed, the canonical product on $SV'$ arises from the canonical coproduct on $SV$ after the fashion familiar from Hopf algebra theory. 

\medbreak 
When $v \in V$ the Fock field operator $\pi (v)$ is defined in terms of the creator $c(v)$ and annihilator $a(v)$ according to the usual prescription $\sqrt{2} \  \pi (v) = c(v) + a(v)$. These operators are initially defined on $SV$ and then extend to $SV'$ by antiduality: thus, if $\Phi \in SV'$ and $\psi \in SV$ then $ [\: c(v) \Phi ] (\psi) = \Phi (a(v) \psi)$ and $ [ \: a(v) \Phi ] (\psi) = \Phi (c(v) \psi)$ so that $ [ \: \pi (v) \Phi ] (\psi) = \Phi (\pi (v) \psi)$. The various operators restrict from $SV'$ to the usual domains in $S[V]$: for instance, $\pi (v)$ restricts from $SV'$ to define an operator that is selfadjoint on the natural domain $\{ \Phi \in S[V] : \pi(v) \Phi \in S[V] \}$; it is not necessary to establish that $\pi (v)$ is essentially selfadjoint on $SV$ and form the unique selfadjoint extension. We remark that in this context, the canonical commutation relations in Heisenberg form hold without qualification on $SV'$ : thus, if $x,y \in V$ then $ [ \pi (x), \pi (y) ] = i \: \IM <x \vert y> I $. 

\medbreak 
The universal implementability of symplectic automorphisms may be established rather directly within this formalism. By definition, a (generalized) Fock implementer for the symplectic automorphism $g$ of $V$ is a (nonzero) linear map $ U: SV \ra SV'$ that intertwines $\pi (v) \in \End SV$ with $\pi (gv) \in \End SV'$ in the sense $ v \in V \Rightarrow U \pi (v) = \pi (gv) U$. It transpires that each symplectic automorphism $g$ of $V$ admits a (generalized) Fock implementer $U$ that is unique up to scalar multiples; moreover $U$ may be recovered from the corresponding displaced vacuum, which is a Gaussian (the exponential of a quadratic) in $SV'$. Of course, if the commutator $[g,i]$ is of Hilbert-Schmidt class then the Gaussian displaced vacuum lies in $S[V]$ and (when scaled appropriately) $U$ determines a unitary operator on $S[V]$ that implements $g$ in the usual sense. 

\medbreak
Of course, the technique of placing a Hilbert space ${\mathbb E}$ between a suitable subspace $E$ and its antidual $E'$ so as to form a triple $E \subset {\mathbb E} \subset E'$ is well established, though the subspace $E$ is typically provided with extra structure (such as that of a nuclear space) and the antidual $E'$ respects this. The case in which $L^2 ({\mathbb R}^n)$ is placed between the Schwartz space ${\mathcal S} ({\mathbb R}^n)$ and the tempered distributions ${\mathcal S}' ({\mathbb R}^n)$ is prototypical, of course. Of more direct relevance to the present paper is work of the Hida group and others on the White Noise Calculus: here, ${\mathbb E}$ is the $L^2$ space of a Gaussian measure on the dual of a nuclear space, $E$ the space of test white noise functionals and $E'$ the space of generalized white noise functionals; see \cite{HY} and \cite{Obata} for detailed accounts. 

\medbreak 
Traditional approaches to the bosonic Fock representation may be found in \cite{BSZ} \cite{Ber} \cite{BraRob}; traditional approaches to the classical Shale theorem may be found in \cite{Araki} \cite{BSZ} \cite{Ber} \cite{BraRob} \cite{Rui} \cite{GBS} \cite{Shale} \cite{Vergne}. The approach taken in these notes, placing bosonic Fock space between the symmetric algebra and its full antidual, is both natural and elegant. The virtues of placing fermionic Fock space between the exterior algebra and its full antidual have already been discussed elsewhere \cite{RIMS}. The task of presenting a similar treatment for Fock spaces over indefinite inner product spaces will be left to a subsequent paper.

\section{Symmetric Fock spaces} 

Let $V$ be a complex Hilbert space with $< \cdot \vert \cdot>$ as its complex inner product and $J = i \cdot $ as its complex structure. Denote by 
	\[SV = \bigoplus_{d \in {\mathbb N}} S^{d}V
\]
\noindent 
its graded symmetric algebra and by $P^{d} : SV \ra S^{d}V$ projection on the summand of homogeneous degree ${d \in {\mathbb N}}$. Recall that $SV$ carries a standard complex inner product $< \cdot \vert \cdot>$ relative to which the homogeneous summands are mutually perpendicular: $1 \in {\mathbb C} = S^{0}V$ is a unit vector and if $x_1, \dots , x_d, y_1, \dots ,y_d \in V$ then 
\begin{equation}
	<x_1 \cdots x_d \vert y_1 \cdots y_d> = \rm Per [<x_a \vert y_b>] = \sum_{\pi} \prod_{k = 1}^{d} <x_k \vert y_{\pi (k)}>
\end{equation}

\noindent 
where $\rm Per$ denotes the permanent of a square matrix and $\pi$ runs over the group comprising all permutations of $1, \dots ,d$. In particular, if $x,y \in V$ then 
\begin{equation}
	<x^d \vert y^d> = d! <x \vert y>^d
\end{equation}

\noindent 
and if $v, x_1, \dots ,x_a, y_1, \dots ,y_b \in V$ then 
\begin{equation}\label{prod}
	\Bigl< \frac{v^{a + b}}{(a + b)!} \vert x_1 \cdots x_a y_1 \cdots y_b \Bigr> = \Bigl<\frac{v^a}{a!} \vert x_1 \cdots x_a \Bigr> \Bigl<\frac{v^b}{b!} \vert y_1 \cdots y_b \Bigr> 
\end{equation} 

\noindent 
whence bilinearity implies that if $v \in V$ and $\phi , \psi \in SV$ then 
\begin{equation}\label{prodabst} 
	\Bigl< \frac{v^{a + b}}{(a + b)!} \vert \phi \psi \Bigr> = \Bigl<\frac{v^a}{a!} \vert \phi \Bigr> \Bigl<\frac{v^b}{b!} \vert \psi \Bigr> .
\end{equation} 

\medbreak 

\begin{theorem} \label{basis} 
If $V$ contains $\{v_1, \dots ,v_m\}$ as a unitary set then $SV$ contains $\{v^D : D \in {\mathbb N}^m \}$ as a unitary set, where if $D = (d_1, \dots ,d_m) \in {\mathbb N}^m$ then 
	\[v^D = \frac{v_1 ^{d_1} \cdots v_m ^{d_m}}{\sqrt{d_1! \cdots d_m!}}.
\]

\end{theorem}
\begin{proof} 
If $A, B \in {\mathbb N}^m$ are distinct then $< v^A \vert v^B> =0$: either $v^A$ or $v^B$ have distinct degrees or each term in the permanent expansion of $< v^A \vert v^B>$ contains a vanishing inner product. If $D = (d_1, \dots ,d_m)$ then the permanent expansion of $< v^D \vert v^D>$ has exactly $d_1! \cdots d_m!$ nonvanishing terms each of which equals $<v_1 \vert v_1>^{d_1} \cdots <v_m \vert v_m>^{d_m}$.
\end{proof} 

\medbreak 
For future reference, we remark that $S^{d}V$ is spanned by the vectors $\{u^d : u \in V \}$: indeed, $S^{d}V$ is certainly spanned by $\{u_1 \cdots u_d : u_1, \dots ,u_d \in V \}$ and polarization yields 
\begin{equation}
	2^{d} d!\: u_1 \cdots u_d = \sum_{\pm \cdots \pm} \pm \cdots \pm (\pm u_1 \cdots \pm u_d)^d.
\end{equation}

\medbreak 
It proves convenient to introduce the set ${\mathcal F} (V)$ comprising all finite-dimensional complex subspaces of $V$ directed by inclusion. Note that $SV$ is the union of its subalgebras $SM$ as $M$ runs over ${\mathcal F} (V)$: 
	\[SV = \bigcup \{SM : M \in {\mathcal F} (V)\}.
\]
\noindent 
When $M \in {\mathcal F} (V)$ we write $P_M : V \ra M$ for orthogonal projection and write 
	\[{\mathcal F}_{M} (V) = \{N \in {\mathcal F} (V) : M \subset N \}. 
\]
\medbreak 

\begin{theorem}
If $M \in {\mathcal F} (V)$ then the functorial extension of $P_M : V \ra M$ is precisely the orthogonal projection $P_M : SV \ra SM$. 
\end{theorem} 
\begin{proof} 
Formulae of the type \eqref{prod} show that if $v_1, \dots ,v_d \in V$ and $z \in M$ then
	\[<z^d \vert (P_M v_1) \cdots (P_M v_d)> = d! \: <z \vert P_M v_1> \cdots <z \vert P_M v_d> = <z^d \vert v_1 \cdots v_d> 
\]
\noindent 
whence the remark following Theorem \ref{basis} shows that $(P_M v_1) \cdots (P_M v_d) - (v_1 \cdots v_d)$ is perpendicular to $S^{d}M$. 

\end{proof} 
\medbreak 
Denote by $SV'$ the full antidual of the symmetric algebra, comprising all antilinear functionals $SV \ra {\mathbb C}$. Note that the standard complex inner product $< \cdot \vert \cdot>$ linearly embeds $SV$ in $SV'$ via the canonical inclusion 
	\[SV \ra SV' : \phi \mapsto < \cdot \vert \phi> .
\]
\medbreak 

\noindent 
When $\Phi \in SV'$ and $d \in {\mathbb N}$ we may consider $\Phi^{d} : \:= \Phi\circ P^{d}$ as an element of either $SV'$ or $(S^{d}V)'$ as convenient. Note that if $\Phi \in SV'$ then
	\[\Phi = \sum_{d \in {\mathbb N}} \Phi^{d} 
\]
\noindent 
for if also $\psi \in SV$ then each sum is actually finite in the following calculation: 
	\[\Phi(\psi) = \Phi \Bigl(\sum_{d \in {\mathbb N}} P^{d} \psi \Bigr) = \sum_{d \in {\mathbb N}} \Phi(P^{d} \psi) = \sum_{d \in {\mathbb N}} \Phi^{d} (\psi). 
\]
\noindent 
Note also that if to each $d \in {\mathbb N}$ is associated an element $\phi^{d} \in S^{d} V$ then the formal series $\sum_{d \in {\mathbb N}} \phi^{d}$ determines an element of $SV'$. 

\medbreak 
Now, let $\Phi \in SV'$. If $M \in {\mathcal F} (V)$ and $d \in {\mathbb N}$ then the finite-dimensionality of $S^{d}M$ guarantees the existence of a unique $\Phi_{M}^{d} \in S^{d}M$ such that $\Phi \vert S^{d}M = < \cdot \vert \Phi_{M}^{d} >$. If also $N \in {\mathcal F}_{M} (V)$ then $P_M \Phi_{N}^{d} = \Phi_{M}^{d}$ for if $\psi \in S^{d}M$ then $P_M \psi = \psi$ and therefore 
	\[< \psi \vert \Phi_{M}^{d} > = \Phi(\psi) = < \psi \vert \Phi_{N}^{d} > = < \psi \vert P_{M} \Phi_{N}^{d}>. 
\]
\noindent 
\medbreak 

In the opposite direction is the following description of the antidual. 

\medbreak 

\begin{theorem} \label{cons}
If to each $M \in {\mathcal F} (V)$ and $d \in {\mathbb N}$ is associated an element $\Phi_{M}^{d} \in S^{d}M$ satisfying the consistency condition 
	\[N \in {\mathcal F}_{M} (V) \Rightarrow P_M \Phi_{N}^{d} = \Phi_{M}^{d}
\]
\noindent 
then there exists a unique $\Phi \in SV'$ such that if $M \in {\mathcal F} (V)$ and $d \in {\mathbb N}$ then 
	\[\Phi \vert S^{d}M = < \cdot \vert \Phi_{M}^{d} >.
\]
\end{theorem} 
\begin{proof} 
For $\psi \in SV$ we define 
	\[\Phi (\psi) = \sum_{d \in {\mathbb N}} <P^{d}\psi \vert \Phi_{M}^{d}> 
\]
\noindent 
where $M \in {\mathcal F} (V)$ is chosen so that $\psi \in SM$. The choice of $M \in {\mathcal F} (V)$ is immaterial: if also $N \in {\mathcal F} (V)$ and $\psi \in SN$ then each of $<P^{d}\psi \vert \Phi_{M}^{d}>$ and $<P^{d}\psi \vert \Phi_{N}^{d}>$ equals $<P^{d}\psi \vert \Phi_{M + N}^{d}>$ by consistency. The rest of the proof is clear. 
\end{proof}

\medbreak 
In fact, the antidual $SV'$ is naturally an algebra. The most elegant way to see this rests on the fact that $SV$ itself is naturally a coalgebra: the diagonal map $V \ra V \oplus V$ induces an algebra homomorphism $SV \ra S(V \oplus V)$ which when followed by the canonical isomorphism $S(V \oplus V) \ra SV \otimes SV$ yields the (cocommutative) coproduct $\Delta: SV \ra SV \otimes SV$. In these terms, the natural (commutative) product in $SV'$ is defined by the rule that if $\Phi, \Psi \in SV'$ and $\theta \in SV$ then 
	\[[\Phi \Psi](\theta) = [\Phi \otimes \Psi](\Delta \theta).
\]
\medbreak 

\begin{theorem} \label{prodcon}
The natural product in $SV'$ is weakly continuous. 
\end{theorem} 

\begin{proof} 
Explicitly, if $(\Phi_\lambda : \lambda \in \Lambda)$ and $(\Psi_\lambda : \lambda \in \Lambda)$ are nets in $SV'$ converging weakly to $\Phi \in SV'$ and $\Psi \in SV'$ respectively then the net $(\Phi_{\lambda} \Psi_{\lambda}  : \lambda \in \Lambda)$ converges to $\Phi \Psi$ in the same sense: if $\theta \in SV$ then 
	\[\lim_{\lambda \in \Lambda} [\Phi_{\lambda} \Psi_{\lambda}](\theta)= [\Phi \Psi](\theta)
\]
\noindent 
for if $\Delta\theta = \sum_{k = 1}^{K} \xi_k \otimes \eta_k $ then 
	\[[\Phi_{\lambda} \Psi_{\lambda}](\theta) = [\Phi_{\lambda} \otimes \Psi_{\lambda}](\Delta \theta) =  \sum_{k = 1}^{K} \Phi_{\lambda}(\xi_k)\Psi_{\lambda}(\eta_k)
\]
\noindent 
which as $\lambda$ runs over $\Lambda$ converges to 
	\[\sum_{k = 1}^{K} \Phi(\xi_k)\Psi(\eta_k) = [\Phi \otimes \Psi](\Delta \theta) = [\Phi \Psi](\theta).
\]
\end{proof} 
\medbreak 
Note that the canonical inclusion $SV \ra SV'$ is an algebra homomorphism: it is enough to see that if $\phi \in S^{a}V$ and $\psi \in S^{b}V$ then $< \cdot \vert \phi \psi > = < \cdot \vert \phi  > < \cdot \vert \psi  >$; this follows from \eqref{prodabst} and the remark immediately after Theorem \ref{basis}. When identified with its image, $SV$ is weakly dense in $SV'$. 

\medbreak 
\begin{theorem} \label{wklim}
$\Phi \in SV'$ is the weak limit in $SV'$ of the net $(\Phi_{M}^{0} + \cdots + \Phi_{M}^{d}: M \in {\mathcal F}(V), d \in {\mathbb N})$ in $SV$. 
\end{theorem} 
\begin{proof} 
For $\psi \in SV$ choose $M_\psi \in {\mathcal F}(V)$ and $d_\psi \in {\mathbb N}$ so that $\psi \in SM_\psi$ and $\psi^d = 0$ when $d > d_\psi$. If $M \in {\mathcal F}(V)$ contains $M_\psi$ and $d \in {\mathbb N}$ exceeds $d_\psi$ then plainly $\Phi(\psi) = <\psi \vert \Phi_{M}^{0} + \cdots + \Phi_{M}^{d}>$. 

\end{proof} 

\medbreak 
The following less elegant formulation of the product in $SV'$ is occasionally useful. 

\medbreak 

\begin{theorem} \label{compprod}
Let $\Phi$ and $\Psi$ lie in $SV'$. If $M \in {\mathcal F}(V)$ and $d \in {\mathbb N}$ then 
	\[[\Phi \Psi]_{M}^{d} = \sum_{a + b = d} \Phi_{M}^{a} \Psi_{M}^{b} .
\]
\end{theorem} 

\begin{proof} 
This follows from \eqref{prodabst} and the remark after Theorem \ref{basis} : if $v \in M$ then as 
	\[\Delta(v^d) = (v \otimes 1 + 1 \otimes v)^d = \sum_{a + b = d} \frac{d!}{a!\:b!}\:v^a \otimes v^b
\]
\noindent 
so 
\begin{eqnarray*}
	\Bigl<\frac{v^d}{d!} \Big\vert [\Phi \Psi ]_{M}^{d} \Bigr> & = & [\Phi \Psi ]\Bigl ( \frac{v^d}{d!}\Bigr ) = (\Phi \otimes \Psi) \Bigl(\sum_{a + b = d} \frac{v^a}{a!} \otimes \frac{v^b}{b!} \Bigr)\\ & = & \sum_{a + b = d} \Phi \Bigl( \frac{v^a}{a!} \Bigr) \Psi \Bigl( \frac{v^b}{b!} \Bigr) = \sum_{a + b = d} \Bigl< \frac{v^a}{a!} \Big\vert \Phi_{M}^{a} \Bigr> \Bigl< \frac{v^b}{b!} \Big\vert \Psi_{M}^{b} \Bigr> \\ & = & \Bigl< \frac{v^d}{d!} \Big\vert \sum_{a + b = d} \Phi_{M}^{a} \Psi_{M}^{b} \Bigr>.   
\end{eqnarray*}

\end{proof} 

\medbreak 
We remark that this actually provides an alternative construction of the product in $SV'$: if $\Phi, \Psi \in SV'$ then the assignment 
	\[M \in {\mathcal F}(V), d \in {\mathbb N} \Rightarrow [\Phi \Psi]_{M}^{d} = \sum_{a + b = d} \Phi_{M}^{a} \Psi_{M}^{b}
\]
\noindent 
is readily confirmed to be consistent in the sense of  \:Theorem \ref{cons}. 
\medbreak 

Now, let $\Phi \in SV'$. To each $M \in {\mathcal F}(V)$ we associate the formal sum 
	\[\Phi_{M} : \:= \sum_{d \in {\mathbb N}} \Phi_{M}^{d} \in SV'
\]
\noindent 
and to this formal sum we associate the number 
	\[\Vert \Phi_{M} \Vert : \:= \sqrt{\Bigl\{\sum_{d \in {\mathbb N}} \Vert \Phi_{M}^{d} \Vert ^2 \Bigr\}} \in [\:0, \infty ].
\]
\noindent 
Let also $N \in {\mathcal F}_{M}(V)$: if $d \in {\mathbb N}$ then the consistency condition $\Phi_{M}^{d} = P_{M} \Phi_{N}^{d}$ implies that $\Vert \Phi_{M}^{d} \Vert \leq \Vert \Phi_{N}^{d} \Vert$ ; consequently, summation yields $\Vert \Phi_{M} \Vert \leq \Vert \Phi_{N} \Vert$. It follows that the net $( \Vert \Phi_{N} \Vert : N \in {\mathcal F}(V) )$ in $[\:0, \infty]$ is increasing, with the same supremum as its subnet $( \Vert \Phi_{N} \Vert : N \in {\mathcal F}_{M}(V) )$ for each $M \in {\mathcal F}(V)$. Define 
	\[ \Vert \Phi \Vert : \: = \sup_{N} \Vert \Phi_{N} \Vert = \lim_{N} \Vert \Phi_{N}\Vert. 
\]
\medbreak 

\begin{theorem} \label{norm}
If $\Phi \in SV'$ then $\Vert \Phi \Vert$ is its operator norm as an antilinear functional on $SV$ in the sense 
	\[\Vert \Phi \Vert = \sup \{ \Vert \Phi (\psi ) \vert : \psi \in SV,\: \Vert \psi \Vert \leq 1 \}.
\]
\end{theorem} 

\begin{proof} 
Let $\psi \in SV$ be a unit vector: if $M \in {\mathcal F}(V)$ and $d \in {\mathbb N}$ are chosen so that $\psi = \psi^0 + \cdots + \psi^d \in SM$ then 
	\[\vert \Phi (\psi) \vert = \Big\vert \Bigl< \psi \Big\vert \sum_{a = 0}^{d} \Phi_{M}^{a} \Bigr> \Big\vert \leq \Big\Vert \sum_{a = 0}^{d} \Phi_{M}^{a} \Big\Vert \leq \Vert \Phi_{M} \Vert 
\]
 
\medbreak 
\noindent 
so the operator norm of $\Phi$ is at most $\Vert \Phi \Vert$. Let $M \in {\mathcal F}(V)$ and $d \in {\mathbb N}$: if $\Phi_{M}^{0} + \cdots + \Phi_{M}^{d}$ is nonzero then the unit vector 
	\[\psi : \:= \Bigl(\sum_{a = 0}^{d} \Phi_{M}^{a} \Bigr) \Big/ \Big\Vert \sum_{a = 0}^{d} \Phi_{M}^{a} \Big\Vert 
\]

\noindent 
satisfies 
	\[\Phi (\psi) = \Bigl< \psi \Big\vert \sum_{a = 0}^{d} \Phi_{M}^{a} \Bigr> = \Big\Vert \sum_{a = 0}^{d} \Phi_{M}^{a} \Big\Vert 
\]

\noindent 
whence the arbitrary nature of $M$ and $d$ implies that the operator norm of $\Phi$ is at least $\Vert \Phi \Vert$.

\end{proof}

 We are now in a position to introduce symmetric Fock space as  
	\[S[V] = \{\Phi \in SV' : \Vert \Phi \Vert < \infty \}.
\]
\medbreak  
\noindent 
Plainly, $S[V]$ is a complex vector space upon which $\Vert \cdot \Vert$ defines a norm. In fact, this norm is induced by a complex inner product: indeed, if $\Phi, \Psi \in S[V]$ and $M \in {\mathcal F}(V)$ then the parallelogram law in homogeneous summands of $SM$ yields 
	\[\Vert (\Phi - \Psi)_{M} \Vert^2 + \Vert (\Phi + \Psi)_{M} \Vert^2 = 2\{\Vert \Phi_{M} \Vert^2 + \Vert \Psi_{M} \Vert^2 \}
\]

\noindent 
whence passage to the supremum as $M$ runs over ${\mathcal F}(V)$ yields 
	\[\Vert \Phi - \Psi \Vert^2 + \Vert \Phi + \Psi \Vert^2 = 2\{ \Vert \Phi \Vert^2 + \Vert  \Psi \Vert^2 \} 
\]

\noindent 
so the parallelogram law holds in $S[V]$. Accordingly, $\Vert \cdot \Vert$ is induced by the inner product $< \cdot \vert \cdot >$ defined by the rule that if $\Phi, \Psi \in S[V]$ then 
	\[<\Phi \vert \Psi> \:= \frac{1}{4} \sum_{p = 0}^{3} i^{-p} \Vert \Phi + i^p \Psi \Vert^2.
\]
\medbreak 
\begin{theorem} \label{pythM}
If $\Phi \in S[V]$ and $M \in {\mathcal F}(V)$ then 
	\[\Vert \Phi \Vert^2 = \Vert \Phi - \Phi_{M} \Vert^2 + \Vert \Phi_{M} \Vert^2. 
\]
\end{theorem} 
\begin{proof} 
Let $N \in {\mathcal F}_{M}(V)$. If $d \in {\mathbb N}$ then consistency and the Pythagorean law in $S^{d}N$ yield 
	\[\Vert \Phi_{N}^{d} \Vert^2 = \Vert \Phi_{N}^{d} - \Phi_{M}^{d}\Vert^2 + \Vert  \Phi_{M}^{d}\Vert^2 = \Vert (\Phi - \Phi_{M})_{N}^{d}\Vert^2 + \Vert ( \Phi_{M})_{N}^{d}\Vert^2 
\]
\medbreak 
\noindent 
whence summation yields 
	\[\Vert \Phi_{N} \Vert^2 = \Vert (\Phi - \Phi_{M})_{N}\Vert^2 + \Vert ( \Phi_{M})_{N}\Vert^2.
\]
\medbreak 
\noindent 
Passage to the supremum as $N$ runs over ${\mathcal F}_{M}(V)$ concludes the argument. 
	
\end{proof} 
\medbreak 
As is readily checked, it is also the case that if $\Phi \in S[V]$ then 
\begin{equation} \label{pars}
	\Vert \Phi \Vert^2 = \sum_{d \in {\mathbb N}} \Vert \Phi^d \Vert^2.	
\end{equation}

\medbreak 
In fact, $SV$ is dense in the inner product space $S[V]$. 
\medbreak 
\begin{theorem}\label{dense}
If $\Phi \in S[V]$ then the net $(\Phi_{M}^{0} + \cdots + \Phi_{M}^{d}: M \in {\mathcal F}(V), d \in {\mathbb N})$ in $SV$ converges to $\Phi$ in $S[V]$. 
\end{theorem} 
\begin{proof} 
Let $\varepsilon > 0$. Choose $M_{\varepsilon} \in {\mathcal F}(V)$ so that $\Vert \Phi_{M_{\varepsilon}} \Vert^2 > \Vert \Phi \Vert^2 - \varepsilon^2$ and choose $d_{\varepsilon} \in {\mathbb N}$ so that $\Vert \Phi_{M_{\varepsilon}}^{0} + \cdots + \Phi_{M_{\varepsilon}}^{d_{\varepsilon}} \Vert^2 > \Vert \Phi \Vert^2 - \varepsilon^2$. If $M \in {\mathcal F}_{M_{\varepsilon}}(V)$ and if $d \in {\mathbb N})$ is at least $d_{\varepsilon} \in {\mathbb N}$ then 
\begin{eqnarray*} 
\Big\Vert \Phi - \sum_{a = 0}^{d} \Phi_{M}^{a} \Big\Vert^2 & = & \sum_{a > d} \Vert \Phi^a \Vert^2 + \sum_{a = 0}^{d} \Vert \Phi^a - \Phi_{M}^{a}\Vert^2 \\ & = & \sum_{a \in {\mathbb N}} \Vert \Phi^a \Vert^2 - \sum_{a = 0}^{d} \Vert  \Phi_{M}^{a}\Vert^2 \\ & = & \Vert \Phi \Vert^2 - \Big\Vert \sum_{a = 0}^{d} \Phi_{M}^{a} \Big\Vert^2 <\: {\varepsilon}^2. 
\end{eqnarray*}
\end{proof} 

\medbreak 
Further, $S[V]$ is actually the Hilbert space completion of $SV$. 

\medbreak 
\begin{theorem} 
The inner product space $S[V]$ is complete. 
\end{theorem} 
\begin{proof} 
Let $( ^{j} \Phi : j \in {\mathbb N})$ be a Cauchy sequence in $S[V]$. If $M \in {\mathcal F}(V)$ and $d \in {\mathbb N}$ then $(^{j}\Phi_{M}^{d} : j \in {\mathbb N})$ is (by domination) a Cauchy sequence in the finite-dimensional (hence complete) space $S^{d}M$ so we may define $\Phi_{M}^{d} : \:= \lim_{j} (^{j}\Phi_{M}^{d})$. If also $N \in {\mathcal F}_{M}(V)$ then $P_{M} ^{j}\Phi_{N}^{d} = \: ^{j}\Phi_{M}^{d}$ so that continuity of $P_{M} : S^{d}N \ra S^{d}M $ implies $P_{M} \Phi_{N}^{d} = \Phi_{M}^{d}$. Now Theorem \ref{cons} furnishes a unique $\Phi \in SV'$ such that if $M \in {\mathcal F}(V)$ and $d \in {\mathbb N}$ then $\Phi \vert S^{d}M = < \cdot \vert \Phi_{M}^{d} >$. Let $\varepsilon > 0$ and choose $j_{\varepsilon} \in {\mathbb N} $ so that if $p, q \geq j_{\varepsilon}$ then $\Vert ^{q}\Phi -  ^{p}\Phi \Vert \leq \varepsilon$. If $M \in {\mathcal F}(V)$ then $\Vert ^{q}\Phi_{M} -  ^{p}\Phi_{M} \Vert \leq \varepsilon$ so that (upon inspection of homogeneous summands) letting $p = j \geq j_{\varepsilon}$ and $q \ra \infty$ results in $\Vert (\Phi -  ^{j}\Phi)_{M} \Vert \leq \varepsilon$; as $M$ is arbitrary, it follows that if $j \geq j_{\varepsilon}$ then $\Vert \Phi -  ^{j}\Phi \Vert \leq \varepsilon$. This places $\Phi$ in $S[V]$ as the limit of $( ^{j} \Phi : j \in {\mathbb N})$. 
\end{proof} 

\medbreak 
Of course, the canonical inclusion $SV \ra S[V]$ is isometric. 
\medbreak 
We shall have occasion to use the following assertion of compatibility. 

\medbreak
\begin{theorem} \label{compat} 
If $\Phi \in S[V]$ and $\psi \in SV$ then $\Phi (\psi) = <\psi \vert \Phi>$. 
\end{theorem} 
\begin{proof} Select $M_{\psi} \in {\mathcal F}(V)$ and $d_{\psi} \in {\mathbb N}$ so that $\psi = \psi^{0} + \cdots + \psi^{d_{\psi}} \in SM_{\psi} $. If $M \in {\mathcal F}_{M_{\psi}}(V)$ and $d \in {\mathbb N}$ exceeds $d_{\psi}$ then $\psi = \psi^{0} + \cdots + \psi^{d_{\psi}} \in SM $ so 
	\[\Phi (\psi) = <\psi \vert \Phi_{M}^{0} + \cdots + \Phi_{M}^{d}>.
\]
\noindent 
An application of Theorem \ref{dense} ends the proof. 
\end{proof} 
\medbreak 
We shall also have need for the subspace of $S[V]$ comprising all elements of homogeneous degree $d \in {\mathbb N}$: 
	\[S^{d}[V] = \{ \Phi \in S[V] : \Phi \circ P^{d} = \Phi \}.
\]

\medbreak
\begin{theorem} \label{closure} 
If $d \in {\mathbb N}$ then $S^{d}[V]$ is precisely the closure of $S^{d}V$ in $S[V]$. 
\end{theorem} 
\begin{proof} 
Plainly, $S^{d}V \subset S^{d}[V]$ and \eqref{pars} implies that the map $S[V] \ra S[V] : \Phi \mapsto \Phi \circ P^{d}$ is continuous, so $\overline{S^{d}V} \subset S^{d}[V]$. For the reverse inclusion, apply Theorem \ref{dense}. 
\end{proof} 

\medbreak 
Note that \eqref{pars} shows that $S[V]$ is the Hilbert space direct sum of its homogeneous subspaces: 
	\[S[V] = \overline{\bigoplus_{d \in {\mathbb N}}} S^{d}[V].
\]
\medbreak 
This notion of homogeneity may be conveniently reformulated and applies to the full antidual $SV'$. Explicitly, the group ${\mathbb R^+}$ of positive reals has a scaling action $\sigma$ on $V$ given by 
	\[t \in {\mathbb R^+} , v \in V \Rightarrow \sigma_{t}v = tv 
\]
\medbreak 
\noindent 
which extends to $SV'$ by functoriality and then to $SV'$ by antiduality, so that 
	\[t \in {\mathbb R^+}, \Phi \in SV', \psi \in SV \Rightarrow [\sigma_{t}\Phi ] (\psi ) = \Phi (\sigma_{t}\psi ). 
\]
\medbreak 
\noindent 
In these terms, the elements of homogeneous degree $d \in {\mathbb N}$ are those on which $\sigma_{t}$ acts as multiplication by $t^{d}$ whenever $t \in {\mathbb R^+}$. 
\medbreak 
In the sequel, our interest will centre largely on quadratics: elements of $S^{2}V'$ . It is convenient to discuss these a little more fully here. Let us say that the (antilinear) map $Z: V \ra V$ is symmetric precisely when 
	\[x,y \in V \Rightarrow <y \vert Zx> = <x \vert Zy>.
\]
\medbreak 
\noindent 
More generally, let us say that the (antilinear) map $Z: V \ra V'$ is symmetric precisely when 
	\[x,y \in V \Rightarrow Zx(y) = Zy(x).
\]
\medbreak
\noindent
Plainly, the space $S^{2}V'$ of all quadratics $\zeta$ is canonically isomorphic to the space of all symmetric antilinear maps $Z: V \ra V'$ via the rule 
\begin{equation} \label{quad}
	x,y \in V \Rightarrow \zeta(xy) = Zx(y). 
\end{equation}

\medbreak 
\begin{theorem} \label{corr}
$S^{2}[V]$ is canonically isomorphic to the space $\Sigma^{2}[V]$ comprising all Hilbert-Schmidt symmetric antilinear maps $V \ra V$: explicitly, $\zeta \in S^{2}[V]$ and $Z \in \Sigma^{2}[V]$ correspond when 
	\[x,y \in V \Rightarrow \zeta(xy) = <y \vert Zx>.
\]
\end{theorem} 
\begin{proof} 
Let $M \in {\mathcal F}(V)$ have $(v_1 , \dots , v_m)$ as unitary basis. Note that $\zeta_{M} \in S^{2}M$ corresponds canonically to the symmetric antilinear map $Z_{M} : M \ra M : v \mapsto (Zv)_{M}$. Accordingly, from Theorem \ref{basis} it follows that 
	\[\zeta_{M} = \frac{1}{2}\sum_{a,b} <v_a v_b \vert \zeta_{M}>v_a v_b = \frac{1}{2}\sum_{a,b} <v_a \vert Z_{M}v_b> v_a v_b 
\]
\medbreak 
\noindent 
whence 
	\[\Vert \zeta_{M} \Vert^2 = \frac{1}{2}\sum_{a,b} \vert<v_a \vert Z_{M}v_b>\vert^2 = \frac{1}{2} \Vert Z_{M} \Vert_{HS}^2.
\]
\medbreak 
\noindent 
Passage to the supremum as $M$ runs over ${\mathcal F}(V)$ now shows not only that $Z$ maps $V$ to itself but also that $Z: V \ra V$ is Hilbert-Schmidt with $\Vert Z \Vert_{HS} = \sqrt{2} \Vert \zeta \Vert$.
\end{proof}

\section{Exponentials, creators and annihilators}
	
Let $v \in V$. We define the creator $c(v): SV \ra SV$ to be the operator of left (equivalently, right) multiplication by $v$: 
	\[\phi \in SV \Rightarrow c(v)\phi = v \phi.
\]
\medbreak 
\noindent 
We define the annihilator $a(v): SV \ra SV$ to be the unique linear derivation such that $a(v) 1 = 0$ and such that if $w \in V$ then $a(v) w = <v \vert w>$. Recall that for $a(v)$ to be a derivation means that if $\phi, \psi \in SV$ then 
	\[a(v) [\phi \psi] = [a(v) \phi] \psi + \phi [a(v) \psi] 
\]

\medbreak
\noindent
so that if $v_1, \dots ,v_m \in V$ then 
	\[a(v) [v_1 \cdots v_m] = \sum_{k = 1}^{m} <v \vert v_k> v_1 \cdots \widehat{v_k} \cdots v_m 
\]
\medbreak
\noindent
where the circumflex $ \:\widehat{\cdot}\: $ signifies omission as usual. 
\medbreak
\begin{theorem} \label{adj}
If $v \in V$ then $c(v)$ and $a(v)$ are mutually adjoint on $SV$ in the sense that if $\phi, \psi \in SV$ then 
	\[<a(v) \phi \vert \psi> = <\phi \vert c(v) \psi>.
\]
\end{theorem}
\begin{proof} 
It is enough to verify the equality when $\phi = x_0 x_1 \cdots x_m$ and $\psi = y_1 \cdots y_m$ for vectors $x_0, x_1, \cdots ,x_m, y_1, \cdots ,y_m$ in $V$; in this case, verification amounts to an elementary permanent expansion.  
\end{proof} 

\medbreak
We extend the definition of creators and annihilators to the antidual $SV'$ by antiduality. Explicitly, let $v \in V$: for $\Phi \in SV'$ and $\psi \in SV$ we define 
	\[[c(v) \Phi](\psi) = \Phi (a(v) \psi) 
\]
	\[[a(v) \Phi](\psi) = \Phi (c(v) \psi).
\]
\medbreak
\noindent
In both the original and this extended context, creators and annihilators satisfy the canonical commutation relations in the following form. 
\medbreak 
\begin{theorem} \label{CCR}
If $x,y \in V$ then 
	\[[a(x), a(y)] = 0
\]
	\[[a(x), c(y)] = <x \vert y> I
\]
	\[[c(x), c(y)] = 0.
\]
\end{theorem}
\begin{proof} 
Validity on $SV'$ follows at once by antiduality from validity on $SV$. Here, the last identity is plain from commutativity of $SV$ while the first then follows by Theorem \ref{adj} ; the central identity holds since if $\phi \in SV$ then 
	\[a(x)c(y)\phi = a(x)[y \:\phi] = [a(x)y]\phi + y[a(x)\phi] = <x \vert y> \phi + c(y)a(x) \phi.
\]

\end{proof} 
\medbreak 
When $SV'$ is given the topology of pointwise convergence, the extended creators and annihilators are continuous. 

\medbreak 
\begin{theorem} \label{wkcon}
If $v \in V$ then $c(v)$ and $a(v)$ are weakly continuous on $SV'$. 
\end{theorem} 
\begin{proof} 
Let $(\Phi_\lambda : \lambda \in \Lambda)$ be a net converging weakly to $\Phi$ in $SV'$: if $\psi \in SV$ then as $\lambda$ runs over $\Lambda$ so 
	\[[c(v)\Phi_\lambda](\psi) = \Phi_\lambda(a(v)\psi) \ra \Phi(a(v)\psi) = [c(v)\Phi](\psi)
\]
\medbreak
\noindent
whence $c(v)\Phi_\lambda \ra c(v)\Phi$ weakly and $a(v)\Phi_\lambda \ra a(v)\Phi$ similarly.
\end{proof} 
\medbreak 
The extended creators and annihilators are indeed extensions of the originals relative to the canonical inclusion $SV \ra SV'$: let $v \in V$ and $\phi \in SV$; if also $\psi \in SV$ then by Theorem \ref{adj} it follows that 
	\[[c(v)<\cdot \vert \phi>](\psi) = <\cdot \vert \phi> (a(v) \psi) = <\cdot \vert c(v) \phi> (\psi) 
\]
\medbreak 
\noindent 
whence $c(v) <\cdot \vert \phi> = <\cdot \vert c(v)\phi>$ and $a(v) <\cdot \vert \phi> = <\cdot \vert a(v)\phi>$ likewise. Further, the extended creators and annihilators inherit the following properties from the originals. 
\medbreak 
\begin{theorem} \label{ca}
If $v \in V$ then $c(v): SV' \ra SV'$ is multiplication by $<\cdot \vert v>$ and $a(v) : SV' \ra SV'$ is a derivation. 
\end{theorem} 
\begin{proof} 
For the annihilator, let $\Phi$ and $\Psi$ lie in $SV'$: Theorem \ref{wklim} furnishes nets $(\phi_\lambda: \lambda \in \Lambda)$ and $(\psi_\lambda: \lambda \in \Lambda)$ in $SV$ converging weakly to $\Phi$ and $\Psi$ respectively, so letting $\lambda$ run over $\Lambda$ in 
	\[a(v) [\phi_\lambda  \psi_\lambda] = [a(v) \phi_\lambda] \psi_\lambda + \phi_\lambda [a(v) \psi_\lambda] 
\]
\medbreak
\noindent
yields the desired equality 
	\[a(v) [\Phi \Psi] = [a(v) \Phi] \Psi + \Phi [a(v) \Psi] 
\]
\medbreak
\noindent 
on account of Theorem \ref{prodcon} and Theorem \ref{wkcon}. For the creator, argue by weak continuity or let $\Phi \in SV'$: if $u \in V$ then 
\begin{eqnarray*}
[<\cdot \vert v> \Phi] (u^d) & = & [<\cdot \vert v> \otimes \Phi] \Bigl(\sum_{a + b = d} \frac{d!}{a!\:b!} u^a \otimes u^b \Bigr) \\ & = & \sum_{a + b = d}\frac{d!}{a!\:b!} <u^a \vert v > \Phi (u^b) \\ & = & d<u \vert v> \Phi (u^{d - 1} ) =  \Phi (d <v \vert u> u^{d - 1} ) \\ & = & \Phi (a(v) u^d)  =  [c(v) \Phi ](u^d) 
\end{eqnarray*}
\medbreak 
\noindent 
whence the discussion after Theorem \ref{basis} implies that $c(v) \Phi = <\cdot \vert v> \Phi$. 

\end{proof} 
\medbreak 
Another familiar inherited property concerns the Fock vacuum $1 \in {\mathbb C} = S^0V$. 
\medbreak 
\begin{theorem} \label{vac}
The antifunctionals in $SV'$ killed by each annihilator are exactly the scalar multiples of $< \cdot \vert 1>$. 
\end{theorem} 
\begin{proof} 
By definition, each annihilator vanishes on $1 \in SV$ and hence on $< \cdot \vert 1> \in SV'$. Conversely, let $\Phi \in SV'$ lie in the kernel of each annihilator. If $v_0, v_1, \dots ,v_m \in V$ then 
	\[0 = [a(v_0) \Phi] (v_1 \cdots v_m) = \Phi (v_0 v_1 \cdots v_m) 
\]
\medbreak 
\noindent 
so that $\Phi$ vanishes on $\oplus_{d > 0} S^d V$ and is therefore proportional to $< \cdot \vert 1>$. 
\end{proof} 

\medbreak 
Similarly or otherwise, it is easily checked that each creator is actually injective. 
\medbreak 
In order to consider creators and annihilators as operators in symmetric Fock space $S[V]$ we investigate their relationship to $\Vert \cdot \Vert$. As preparation, let $v \in V$ and let 
	\[\Phi = \sum_{d \in {\mathbb N}}\Phi^d \in SV'.
\]
\medbreak 
\noindent 
Plainly, if $d \in {\mathbb N}$ then $(c(v)\Phi)^{d + 1} = c(v) \Phi^d$ and $(a(v)\Phi)^d = a(v) \Phi^{d + 1}$. Let also $M \in {\mathcal F} (V)$ contain $v$: if $\psi \in S^{d + 1}V$ then 
\begin{eqnarray*}
	(c(v) \Phi)^{d + 1} (\psi) & = & c(v) \Phi^d (\psi) = \Phi^d (a(v) \psi)\\ & = & <a(v) \psi \vert \Phi_{M}^d > = <\psi \vert c(v)\Phi_{M}^d> 
\end{eqnarray*}

\medbreak 
\noindent 
whence 
	\[(c(v)\Phi)_{M}^{d + 1} = c(v)\Phi_{M}^d 
\]
\medbreak 
\noindent 
and similarly 
	\[(a(v)\Phi)_{M}^d = a(v)\Phi_{M}^{d + 1}.
\]
\medbreak 
\begin{theorem} \label{pyth}
If $v \in V$ and $\Phi \in SV'$ then 
	\[\Vert c(v)\Phi \Vert^2 = \Vert a(v)\Phi \Vert^2 + \Vert v \Vert^2 \Vert \Phi \Vert^2. 
\]
\end{theorem} 
\begin{proof} 
Of course, both sides of the putative equality are numbers in $[0, \infty]$. If $M \in {\mathcal F} (V)$ then Theorem \ref{adj} and the canonical commutation relations in Theorem \ref{CCR} imply that 	
\begin{eqnarray*}
	\Vert c(v) \Phi_{M}^d \Vert^2 & = & <c(v) \Phi_{M}^d \vert c(v) \Phi_{M}^d> = <\Phi_{M}^d \vert a(v) c(v) \Phi_{M}^d > \\ & = & < \Phi_{M}^d \vert c(v) a(v) \Phi_{M}^d + <v \vert v> \Phi_{M}^d > \\ & = & \Vert a(v) \Phi_{M}^d \Vert^2 + \Vert v \Vert^2 \Vert \Phi_{M}^d \Vert^2 
\end{eqnarray*}
\medbreak 
\noindent 
whence the formulae derived prior to the Theorem imply that if $M$ contains $v$ then  
	\[\Vert (c(v)\Phi)_{M}^{d + 1} \Vert^2 = \Vert (a(v)\Phi)_{M}^{d - 1} \Vert^2 + \Vert v \Vert^2 \Vert \Phi_{M}^d \Vert^2.
\]
\medbreak
\noindent 
Summation over $d > 0$ together with the evident equalities $\Vert (c(v)\Phi)_{M}^1 \Vert = \Vert c(v)\Phi_{M}^0 \Vert = \Vert v \Vert \Vert \Phi_{M}^0 \Vert$ and $(c(v)\Phi)_{M}^0 = 0$ yields 
	\[\Vert (c(v)\Phi)_{M} \Vert^2 = \Vert (a(v)\Phi)_{M} \Vert^2 + \Vert v \Vert^2 \Vert \Phi_{M} \Vert^2. 
\]
\medbreak 
\noindent 
Passage to the supremum as $M$ runs over ${\mathcal F} (V)$ while containing $v$ concludes the proof. 
\end{proof}
\medbreak 
When $v \in V$ we may now consider $c(v)$ and $a(v)$ as operators in $S[V]$: thus, $c(v)$ has natural domain $\{ \Phi \in S[V]: c(v)\Phi \in S[V] \}$ and  $a(v)$ has natural domain $\{ \Phi \in S[V]: a(v)\Phi \in S[V] \}$. Note that these domains coincide by Theorem \ref{pyth} and plainly contain $SV$. 
\medbreak 
\begin{theorem} \label{adjoints}
When $v \in V$ the operators $c(v)$ and $a(v)$ in $S[V]$ are mutual adjoints: $c(v)^* = a(v)$ and $a(v)^* = c(v)$. 
\end{theorem} 
\begin{proof} 
To see that $a(v) \subset c(v)^*$ let $\Phi$ and $\Psi$ lie in the domain of $a(v)$ and $c(v)$. If $d > 0$ and $M \in {\mathcal F} (V)$ contains $v$ then 
\begin{eqnarray*}
	\Bigl< \sum_{a = 0}^{d - 1} \Psi_{M}^a \Big\vert a(v)\Phi \Bigr> & = & [a(v)\Phi]\Bigl(\sum_{a = 0}^{d - 1} \Psi_{M}^a \Bigr) \\ & = & \Phi \Bigl( c(v)\sum_{a = 0}^{d - 1} \Psi_{M}^a \Bigr) \\ & = & \Phi \Bigl( \sum_{a = 0}^{d} (c(v)\Psi)_{M}^a \Bigr) \\ & = & \Bigl<\sum_{a = 0}^{d} (c(v)\Psi)_{M}^a \Big\vert\Phi \Bigr> 
\end{eqnarray*}
\medbreak
\noindent
whence it follows by Theorem \ref{dense} that 
	\[< \Psi \vert a(v)\Phi> = <c(v)\Psi \vert \Phi>.
\]
\medbreak 
\noindent 
To see that $c(v)^* \subset a(v)$ let $\Phi$ lie in the domain of $c(v)^*$. If $\psi \in SV$ then Theorem \ref{compat} implies that 
\begin{eqnarray*}
	[c(v)^* \Phi](\psi) & = & <\psi \vert c(v)^* \Phi> = <c(v)\psi \vert \Phi>\\ & = & \Phi (c(v)\psi) = [a(v)\Phi](\psi)
\end{eqnarray*}
\medbreak 
\noindent 
whence $a(v)\Phi = c(v)^*\Phi \in S[V]$. Thus $c(v)^* = a(v)$; likewise $a(v)^* = c(v)$.
\end{proof} 
\medbreak 
As a corollary, the operators $c(v)$ and $a(v)$ in $S[V]$ are closed: more directly, this may be seen as follows. Let $(\Phi_j : j \in {\mathbb N})$ be a sequence in the domain of $c(v)$ such that as $j \ra \infty$ both $\Phi_j \ra \Phi$ and $c(v)\Phi_j \ra \Psi$ in $S[V]$. On the one hand, as $\Phi_j \ra \Phi$ in $S[V]$ so $\Phi_j \ra \Phi$ in $SV'$ by Theorem \ref{compat} and therefore $c(v)\Phi_j \ra c(v)\Phi$ in $SV'$ by Theorem \ref{wkcon}; on the other hand, $c(v)\Phi_j \ra \Psi$ in $S[V]$ and hence in $SV'$. Thus $c(v)\Phi = \Psi \in S[V]$ and so $c(v)$ is closed. 
\medbreak 
We shall require certain precise estimates for the norms of a creator and its powers on homogeneous elements of $S[V]$. For these, let $v \in V$ be (without loss) a unit vector. Let $\phi \in S^d V$ and choose $M \in {\mathcal F} (V)$ so that $v \in M$ and $\phi \in S^d M$. Extend $v = v_0$ to a unitary basis $(v_0, v_1, \dots ,v_m)$ for $M$. From Theorem \ref{basis} it follows that 
	\[\phi = \sum_{D} \phi_{D} \frac{v_{0}^{d_0} v_{1}^{d_1} \cdots v_{m}^{d_m}}{\sqrt{d_0 ! d_1 ! \cdots d_m !}}
\]
\medbreak 
\noindent 
where summation extends over all multiindices $D = (d_0, d_1, \dots ,d_m)$ with $d_0 + d_1 + \cdots + d_m = d$. Now 
		\[v \phi = \sum_{D} \sqrt{ \frac{(d_0 + 1)!}{d_0 !}} \: \phi_{D} \frac{v_{0}^{d_0 + 1} v_{1}^{d_1} \cdots v_{m}^{d_m}}{\sqrt{(d_0 + 1) ! d_1 ! \cdots d_m !}}
\]
\noindent 
so 
	\[\Vert v \phi \Vert^2 = \sum_{D} (d_0 + 1) \vert \phi_{D} \vert^2 \leq (d + 1) \Vert \phi \Vert^2. 
\]
\medbreak 
\noindent
This elementary estimate is the basis for the following result. 

\medbreak 
\begin{theorem} \label{binom}
Let $v \in V$ and let $a \in {\mathbb N}$. If $b \in {\mathbb N}$ and $\Phi \in S^b [V]$ then 
	\[\Vert v^a \Phi \Vert^2 \leq \frac{(a + b)!}{a! \: b!} \Vert v^a \Vert^2 \Vert \Phi \Vert^2. 
\]
\end{theorem} 
\begin{proof} 
Allowing $v$ to have arbitrary norm, the inequality immediately prior to the Theorem shows that if $\phi \in S^b V$ then 
	\[\Vert v \phi \Vert^2 \leq (b + 1) \Vert v \Vert^2 \Vert \phi \Vert^2 
\]
\medbreak 
\noindent 
whence induction shows that 
	\[\Vert v^a \phi \Vert^2 \leq (a + b) \cdots (1 + b) \Vert v \Vert^{2 a} \Vert \phi \Vert^2 = \frac{(a + b)!}{a! \: b!} \Vert v^a \Vert^2 \Vert \phi \Vert^2.
\]
\medbreak 
\noindent 
Thus, if $\Phi \in S^b [V]$ and $M \in {\mathcal F} (V)$ contains $v$ then 
	\[\Vert (v^a \Phi)_{M} \Vert^2 \leq \frac{(a + b)!}{a! \: b!} \Vert v^a \Vert^2 \Vert \phi_{M} \Vert^2
\]
\medbreak 
\noindent 
and so passage to the supremum confirms the claimed equality. 
\end{proof} 
\medbreak 
In fact, if $v \in V$ and $a, b \in {\mathbb N}$ then the operator norm of $c(v)^a : S^b [V] \ra S^{a + b} [V]$ is exactly $\sqrt{(a + b)!/a! \: b!} \Vert v^a \Vert $ as may be checked by computing $\Vert c(v)^a v^b \Vert$. 
\medbreak 
Regarding exponentials let us begin simply, considering first the exponentials in $SV'$ of vectors in $V$. To be precise, when $z \in V$ we define 
\begin{equation}
	e^z : \: = \sum_{n \in {\mathbb N}} \frac{z^n}{n!} \in SV'. 
\end{equation}
\medbreak 
\noindent 
As usual, this formal power series is (in the first instance) weakly convergent, for individual elements of $SV$ vanish in sufficiently high degrees. 
\medbreak 
These simple exponentials are called coherent vectors; they are common eigenvectors for the annihilators. 
\medbreak 
\begin{theorem} 
If $v$ and $z$ lie in $V$ then 
	\[a(v) [e^z] = <v \vert z> e^z .
\]
\end{theorem} 
\begin{proof} 
As $a(v)$ is a derivation, if $n \in {\mathbb N}$ then $a(v) [z^n] = n  <v \vert z> z^{n - 1}$ so 
	\[a(v) \Bigl[ \frac{z^n}{n!} \Bigr] = <v \vert z> \frac{z^{n - 1}}{(n - 1)!}
\]
\medbreak 
\noindent 
from which the Theorem follows upon summation by virtue of the weak continuity expressed in Theorem \ref{wkcon}. 
\end{proof} 

\medbreak 
In fact, these simple exponentials converge not only in $SV'$ but also in $S[V]$. 

\medbreak 
\begin{theorem} \label{enorm}
If $z \in V$ then the coherent vector  $e^z$ lies in $S[V]$ and 
	\[\Vert e^z \Vert^2 = e^{\Vert z \Vert^2} .
\]
\end{theorem} 
\begin{proof} 
If $M \in {\mathcal F} (V)$ contains $z$ then of course $(e^z)_M = e^z$ and 
	\[\Vert (e^z)_M  \Vert^2 = \sum_{n \in {\mathbb N}} \frac{\Vert z^n \Vert^2}{(n!)^2} = \sum_{n \in {\mathbb N}} \frac{\Vert z \Vert^{2 n}}{n!} = e^{\Vert z \Vert^2} .
\]
\medbreak 
\noindent 
Now pass to the supremum as $M$ runs over ${\mathcal F} (V)$ while containing $z$. 

\end{proof} 
\medbreak 
More generally, if $x,y \in V$ then the coherent vectors $e^x$ and $e^y$ have inner product
\begin{equation}
	<e^x \vert e^y> = e^{<x \vert y>} .
\end{equation}
\medbreak 

\begin{theorem} \label{total}
The coherent vectors $\{ e^z : z \in V \}$ constitute a linearly independent total set in $S[V]$.  	
\end{theorem} 
\begin{proof} 
Let $z_1, \dots ,z_m \in V$ be distinct and assume that $\lambda_1, \dots ,\lambda_m \in {\mathbb C}$ are such that 
	\[\lambda_1 e^{z_1} + \cdots + \lambda_m e^{z_m} = 0 
\]
\medbreak 
\noindent 
whence the taking of homogeneous components yields  
	\[d \in {\mathbb N} \Rightarrow \lambda_1 z_{1}^d + \cdots + \lambda_m z_{m}^d = 0.
\]
\medbreak 
\noindent 
Select $v \in V$ outside the finite union 
	\[\bigcup \{ \ker < \cdot \vert z_q - z_p > : 1 \leq p < q \leq m \} 
\]
\medbreak 
\noindent 
of hyperplanes, so that the complex numbers $<v \vert z_1> , \dots , <v \vert z_m>$ are distinct. Now 
	\[d \in {\mathbb N} \Rightarrow <v^d \vert \lambda_1 z_{1}^d + \cdots + \lambda_m z_{m}^d > = 0
\]
\medbreak 
\noindent 
so 
	\[d \in {\mathbb N} \Rightarrow <v \vert z_1>^{d}\lambda_1 + \cdots + <v \vert z_m>^{d}\lambda_m = 0. 
\]
\medbreak 
\noindent 
This Vandermonde system forces the vanishing of $\lambda_1, \dots ,\lambda_m$. This proves that the coherent vectors are linearly independent; we prove that their linear span is dense in $S[V]$ as follows. Let 
	\[\Phi = \sum_{d \in {\mathbb N}} \Phi^d \in S[V] 
\]
\medbreak 
\noindent 
and suppose that $<\Phi \vert e^z > = 0$ whenever $z \in V$. If $z \in V$ is fixed and $\lambda \in {\mathbb C}$ varies then 
	\[0 = < \Phi \vert e^{\lambda z}> = \sum_{d \in {\mathbb N}} <\Phi^d \vert z^d > \frac{\lambda^d}{d!}	
\]
\medbreak \noindent 
whence equating coefficients shows that if $d \in {\mathbb N}$ then $<\Phi^d \vert z^d > = 0$ . As $z \in V$ is arbitrary, so $\Phi$ vanishes in each degree, on account of Theorem \ref{closure} and the discussion following Theorem \ref{basis}. 
\end{proof}

\medbreak 
Of special importance are Gaussians: the exponentials of quadratics. Let the quadratic $\zeta \in S^{2}V'$ correspond to the symmetric antilinear map $Z: V \ra V'$ according to \eqref{quad}. We define the associated Gaussian by 
\begin{equation}
	e^Z = \exp(\zeta) = \sum_{n \in {\mathbb N}} \frac{\zeta^n}{n!} \in SV'
\end{equation}
\medbreak 
\noindent 
where the formal series converges weakly because individual elements of $SV$ vanish in sufficiently high degree. 
\medbreak 
On Gaussians, annihilators act essentially as creators. 
\medbreak 

\begin{theorem} \label{c=a}
If $v \in V$ and if $Z: V \ra V'$ is symmetric antilinear then 
	\[a(v) e^Z = (Zv) e^Z.
\]
\end{theorem} 
\begin{proof} 
Let $Z$ correspond to the quadratic $\zeta \in S^{2}V'$ as usual: the rule \eqref{quad} implies that $a(v)\zeta = Zv$; hence Theorem \ref{ca} implies that if $n \in {\mathbb N}$ then $a(v)\zeta^n = n(Zv)\zeta^{n - 1}$ so Theorem \ref{prodcon} and Theorem \ref{wkcon} imply that $a(v) \exp \zeta = (Zv) \exp \zeta$ .
\end{proof} 
\medbreak 
Contrary to the case for coherent vectors, Gaussians do not automatically lie in symmetric Fock space: in fact, we claim that $e^Z$ lies in $S[V]$ precisely when $Z$ is of Hilbert-Schmidt class and has operator norm strictly less than unity. 
\medbreak 
In order to establish this claim, it is convenient to begin by supposing that $V$ is finite-dimensional. In this case, let $Z: V \ra Z$ be a symmetric antilinear map and note that $Z^2$ is then a selfadjoint (indeed, positive) complex-linear map: 
	\[v \in V \Rightarrow <v \vert Z^2 v> = \Vert Zv \Vert^2. 
\]
\medbreak
\noindent
By diagonalization, $V$ has a unitary basis $(v_1, \dots ,v_m)$ such that if $1 \leq k \leq m$ then $Zv_k = \lambda_k v_k$ with $\lambda_k \geq 0$; in these terms, 
	\[{\rm Det} (I - Z^2) = (1 - {\lambda_1}^2) \cdots (1 - {\lambda_m}^2)
\]
	\[\Vert Z \Vert = \max (\lambda_1 , \dots \lambda_m). 
\]
\medbreak 
\noindent 
The quadratic $\zeta \in S^2 V$ to which $Z$ corresponds canonically is given by 
	\[\zeta = \frac{1}{2} \sum_{k = 1}^{m} \lambda_k v_{k}^2
\]
\medbreak
\noindent
so that if $n \in {\mathbb N}$ then 
	\[\zeta^n = \sum_{N} \binom{n}{n_1 \cdots n_m} \Bigl(\frac{\lambda_1}{2}\Bigr)^{n_1} \cdots \Bigl(\frac{\lambda_m}{2}\Bigr)^{n_m} \: v_{1}^{2 n_1} \cdots v_{m}^{2 n_m} 
\]
\medbreak 
\noindent 
and 
	\[\frac{\Vert \zeta^n \Vert^2}{(n!)^2} = \sum_{N} \binom{2 n_1}{n_1} \cdots \binom{2 n_m}{n_m} \Bigl(\frac{\lambda_1}{2}\Bigr)^{2 n_1} \cdots \Bigl(\frac{\lambda_m}{2}\Bigr)^{2 n_m}
\]
\medbreak 
\noindent 
where summation takes place over all multiindices $N = (n_1 , \dots , n_m) \in {\mathbb N}^m$ for which $n = n_1 + \cdots + n_m$. Consequently, 
\begin{eqnarray*}
	\Vert \exp \zeta \Vert^2 & = & \sum_{n \in {\mathbb N}} \frac{\Vert \zeta^n \Vert^2}{(n!)^2} \\ & = & \sum_{n_1 \in {\mathbb N}} \binom{2 n_1}{n_1} \Bigl(\frac{\lambda_1}{2}\Bigr)^{2 n_1} \cdots \sum_{n_m \in {\mathbb N}} \binom{2 n_m}{n_m} \Bigl(\frac{\lambda_m}{2}\Bigr)^{2 n_m} \\ & = & (1 - \lambda_{1}^2)^{-\frac{1}{2}} \cdots (1 - \lambda_{m}^2)^{-\frac{1}{2}} \\ & = & {\rm Det}^{\frac{1}{2}} (I - Z^2)^{-1} 
\end{eqnarray*}
\medbreak 
\noindent 
provided that each of the nonnegative numbers $\lambda_1 , \dots ,\lambda_m$ is strictly less than unity. 

\medbreak 
We may now establish the claim in full generality. 
\medbreak 
\begin{theorem}\label{Gauss}
If $Z \in \Sigma^2 [V]$ and $\Vert Z \Vert < 1$ then $e^Z \in S[V]$ and 
	\[\Vert e^Z \Vert^2 = {\rm Det}^{\frac{1}{2}} (I - Z^2)^{-1}.
\]

\end{theorem}
\begin{proof}
Let $\zeta \in S^2 [V]$ be the canonical correspondent to $Z \in \Sigma^2 [V]$ as in \eqref{quad}. If $M \in {\mathcal F}(V)$ and if $\zeta_{M} \in S^2 M$ corresponds to $Z_M : M \ra M$ then (by Theorem \ref{compprod} say) $(\exp \zeta)_M = \exp(\zeta_M)$ so that $(e^Z)_M = e^{Z_M}$ while $\Vert Z_M \Vert \leq \Vert Z \Vert < 1$. The finite-dimensional calculation prior to the Theorem yields 
	\[\Vert (e^Z)_M \Vert^2 = {\rm Det}^{\frac{1}{2}} (I - {Z_M}^2)^{-1}.
\]
\medbreak
\noindent
On the one hand, the net $(\Vert (e^Z)_M \Vert : M \in {\mathcal F}(V))$ is increasing by the discussion prior to Theorem \ref{norm} and indeed converges to $\Vert e^Z \Vert$ by definition; on the other hand, the limit of the net $({\rm Det} (I - {Z_M}^2) : M \in {\mathcal F}(V))$ is ${\rm Det} (I - Z^2)$ by trace-norm continuity (or very definition) of the Fredholm determinant. 
\end{proof} 
\medbreak 
Conversely, let $Z: V \ra V'$ be symmetric antilinear with correspondent $\zeta \in S^2 V'$ and suppose that $e^Z \in S[V]$. If $M \in {\mathcal F}(V)$ then the proof of Theorem \ref{corr} yields 
	\[\Vert Z_M \Vert_{HS}^2 = 2 \Vert \zeta_M \Vert^2 < 2 \Vert \exp \zeta_M \Vert^2 = 2 \Vert e^{Z_M} \Vert^2 \leq 2 \Vert e^Z \Vert^2
\]
\medbreak
\noindent
whence $Z \in \Sigma^2 [V]$. Further, $\Vert Z \Vert < 1$: if $\Vert Z \Vert \geq 1$ then let $u \in V$ be an eigenvector for $Z$ with eigenvalue $\lambda \geq 1$; setting $M = {\mathbb C} u \in {\mathcal F}(V)$ yields 
	\[\Vert (e^Z)_M \Vert^2 = \sum_{n \in {\mathbb N}} \binom{2 n}{n} \Bigl(\frac{\lambda}{2}\Bigr)^{2 n} = \infty 
\]
\medbreak
\noindent
which places $e^Z$ outside $S[V]$. 
\medbreak
More generally, we may explicitly compute the inner product between a pair of Gaussians in symmetric Fock space as follows. 
\medbreak 
\begin{theorem} \label{pair}
If $X$ and $Y$ in $\Sigma^2 [V]$ have operator norms strictly less than unity then 
	\[< e^X \vert e^Y > = {\rm Det}^{\frac{1}{2}} (I - YX)^{-1}.
\]
\end{theorem}
\begin{proof}
It follows from Theorem \ref{Gauss} by the principle of analytic continuation that if $M \in {\mathcal F}(V)$ then 
	\[< e^{X_M} \vert e^{Y_M} > = {\rm Det}^{\frac{1}{2}} (I - Y_M X_M)^{-1}
\]
\medbreak
\noindent
since both sides are respectively (antiholomorphic, holomorphic) in $(X_M , Y_M)$ and agree when $X_M = Y_M$. Now pass to the limit as $M$ runs over ${\mathcal F}(V)$ while taking into account Theorem \ref{pythM} and continuity of the determinant. 
\end{proof}
\medbreak
Incidentally, it is perhaps worth recording a related formula. Let $z \in V$ and let $\zeta \in S^2 [V]$ correspond to $Z \in \Sigma^2 [V]$ with $\Vert Z \Vert < 1$. By induction, if $n \in {\mathbb N}$ then 
	\[<z^{2 n} \vert \zeta^n> = (2 n)! \: \Bigl( \frac{1}{2} <z \vert Zz > \Bigr)^n 
\]
\medbreak
\noindent
so that by summation 
\begin{equation}
	< e^z \vert e^Z > = \exp \Bigl( \frac{1}{2} <z \vert Zz > \Bigr).
\end{equation}
\medbreak
\begin{theorem} \label{smooth}
Let $Z \in \Sigma^2 [V]$ and let $\Vert Z \Vert < 1$. If $\phi \in SV$ then $\phi \: e^Z \in S[V]$. 
\end{theorem} 

\begin{proof} 
As usual, linearity and polarization grant us the right to suppose that $\phi = v^n$ for $v \in V$ a unit vector and $n \in {\mathbb N}$. Let $Z$ correspond to $\zeta \in S^2 [V]$ and choose $s > 1$ so that $\Vert sZ \Vert < 1$. From Theorem \ref{binom} it follows at once that 
	\[\Vert v^n e^Z \Vert^2 \leq \sum_{k \in {\mathbb N}}(2k + n) \cdots (2k + 1) \frac{\Vert \zeta^k \Vert^2}{(k!)^2}.
\]
\medbreak
\noindent
Now, the power series 
	\[\sum_{k = 0}^{\infty} \frac{\Vert \zeta^k \Vert^2}{(k!)^2}\: t^{2k}
\]
\medbreak
\noindent 
and 
	\[\sum_{k = 0}^{\infty} (2k + n) \cdots (2k + 1) \frac{\Vert \zeta^k \Vert^2}{(k!)^2}\: t^{2k}
\]
\medbreak
\noindent
have the same radius of convergence; the former converges when $t = s > 1$ so the latter necessarily converges at $t = 1$. 
\end{proof}
\medbreak
We can say a little more about $e^Z$ when $Z \in \Sigma^2 [V]$ and $\Vert Z \Vert < 1$: from Theorem \ref{smooth} it follows that $e^Z$ lies in the domain of each creator polynomial; in fact, by Theorem \ref{c=a} it follows further that $e^Z$ lies in the domain of each polynomial in creators and annihilators. 
\medbreak
So far as symmetric Fock space itself is concerned, there is little point to considering the exponentials of homogeneous elements in $SV'$ having degree greater than two: such exponentials do lie in $SV'$ of course, but they only lie in $S[V]$ when the homogeneous element is zero. 
\medbreak
\begin{theorem}
Let $\zeta \in S^d V'$ be homogeneous of degree $d > 2$. If $\exp \zeta$ lies in $S[V]$ then $\zeta = 0$. 
\end{theorem}
\begin{proof}
If $\exp \zeta$ lies in $S[V]$ then of course its degree $d$ component $\zeta$ lies in $S^d [V]$. Let $v \in V$ be a unit vector and let $M = {\mathbb C} v \in {\mathcal F}(V)$ so that $\zeta_M = \lambda v^d$ for some $\lambda \in {\mathbb C}$: from 
	\[\Vert \exp(\lambda v^d) \Vert^2 = \sum_{n \in {\mathbb N}} \frac{\Vert(\lambda v^d)^n \Vert^2}{(n!)^2} = \sum_{n \in {\mathbb N}} \frac{(dn)!}{(n!)^2}\:\vert \lambda \vert^{2n}
\]
\medbreak
\noindent
and 
	\[\Vert \exp \zeta_M \Vert \leq \Vert \exp \zeta \Vert < \infty
\]
\medbreak
\noindent
it follows that $\lambda = 0$ whence 
	\[<v^d \vert \zeta> = <v^d \vert \zeta_M> = <v^d \vert \lambda v^d> = d! \:\lambda = 0.	
\]
\medbreak
\noindent
To complete the proof, invoke Theorem \ref{closure} in conjunction with the remark following Theorem \ref{basis}. 

\end{proof}

\section{Generalized Fock implementation}

The imaginary part $\Omega$ of the complex inner product $< \cdot \vert \cdot>$ on $V$ is a real symplectic form: an alternating real-bilinear form that is (strongly) nonsingular in the sense that the correspondence $v \leftrightarrow \Omega(v, \cdot)$ is an isomorphism between $V$ and its real dual. The corresponding symplectic group ${\rm Sp}(V)$ comprises all real-linear automorphisms $g$ of $V$ that are symplectic in the sense 
	\[x,y \in V \Rightarrow \Omega(gx,gy) = \Omega(x,y).
\]
\medbreak
\noindent
Note that each $g \in {\rm Sp}(V)$ is automatically bounded: as may be verified by direct calculation, its adjoint relative to the real inner product $(\cdot \vert \cdot) = {\rm Re} < \cdot \vert \cdot >$ on $V$ is given by $g^* = -J g^{-1} J$. 
\medbreak
As is the case for any real-linear endomorphism of a complex vector space, each $g \in {\rm Sp}(V)$ decomposes uniquely as $g = C_g + A_g$ where $C_g = \frac{1}{2}(g - JgJ)$ is complex-linear and $A_g = \frac{1}{2}(g + JgJ)$ is antilinear. 
\medbreak 
\begin{theorem} \label{CA*}
If $g \in {\rm Sp}(V)$ then ${C_g}^* = C_{g^{-1}}$ and ${A_g}^* = -A_{g^{-1}}$ where adjunction is relative to the real inner product $(\cdot \vert \cdot)$ on $V$. 
\end{theorem} 
\begin{proof} 
This follows at once from the formulae for $C_g$ and $A_g$ displayed prior to the Theorem, since $J$ is skew-adjoint and $g^* = -J g^{-1} J$.
\end{proof} 
\medbreak
In terms of the complex inner product $< \cdot \vert \cdot>$ itself, if $g \in {\rm Sp}(V)$ and $x,y \in V$ then 
	\[<C_g x \vert y>\: = \:<x \vert C_{g^{-1}}y> 
\]
	\[<x \vert A_g y> + <y \vert A_{g^{-1}}x> \:= 0.
\]
\medbreak 
\begin{theorem} \label{CA}
If $g \in {\rm Sp}(V)$ then 
	\[C_{g^{-1}} C_g + A_{g^{-1}} A_g = I
\]
	\[A_{g^{-1}} C_g + C_{g^{-1}} A_g = O.
\]
\end{theorem}
\begin{proof}
This is actually valid for any real-linear automorphism $g$ of $V$ and follows upon taking complex-linear and antilinear parts in 
	\[(C_{g^{-1}} + A_{g^{-1}}) (C_g + A_g) = g^{-1} g = I. 
\]
\end{proof}
\medbreak
Note that Theorem \ref{CA*} and Theorem \ref{CA} together imply that if $g \in {\rm Sp}(V)$ then 
	\[{C_g}^* C_g - {A_g}^* A_g = I 
\]
	\[{A_g}^* C_g = {C_g}^* A_g.
\]
\medbreak
\noindent
Thus, ${C_g}^* A_g$ is real self-adjoint and if $v \in V$ then 
	\[\Vert C_g v \Vert^2 = \Vert A_g v \Vert^2 + \Vert v \Vert^2. 
\]
\medbreak
\begin{theorem} \label{Cinv}
If $g \in {\rm Sp}(V)$ then its complex-linear part $C_g$ is invertible. 
\end{theorem} 
\begin{proof} 
The formula immediately prior to the Theorem shows that $C_g$ is injective and indeed bounded below by unity. Similarly $C_{g^{-1}}$ is injective, so Theorem \ref{CA*} implies that $C_g$ has dense range. Together, these facts force $C_g$ to be invertible. 
\end{proof} 
\medbreak
This justifies associating to each $g \in {\rm Sp}(V)$ the antilinear operator 
\begin{equation} \label{Z}
Z_g = -A_g C_g^{-1} = C_{g^{-1}}^{-1} A_{g^{-1}}
\end{equation}
\medbreak
\noindent
which is symmetric antilinear and has operator norm strictly less than unity by virtue of the formulae recorded after Theorem \ref{CA}. 
\medbreak
We shall find it convenient to introduce transformed creators and annihilators. Thus, let $g \in {\rm Sp}(V)$: for $v \in V$ we define 
 \[	c_g(v) = c(C_g v) + a(A_g v) 
\]
 \[	a_g(v) = a(C_g v) + c(A_g v) 
\]
\medbreak
\noindent
as operators on $SV$ and $SV'$. These transformed creators and annihilators continue to satisfy the canonical commutation relations. 
\medbreak
\begin{theorem} \label{CCR'} 
If $g \in {\rm Sp}(V)$ and $x,y \in V$ then 
	\[[a_g (x), a_g (y)] = 0 
\]
	\[[a_g (x), c_g (y)] = <x \vert y> I
\]
	\[[c_g (x), c_g (y)] = 0.
\]
\end{theorem}
\begin{proof} 
Simple application of Theorem \ref{CA*} and Theorem \ref{CA} to the canonical commutation relations of Theorem \ref{CCR}: taking the central identity for example, 
\begin{eqnarray*}
	[a_g (x), c_g (y)] & = & [a(C_g x), c(C_g y)] + [c(A_g x), a(A_g y)] \\ & = & \{ <C_g x \vert C_g y> - <A_g y \vert A_g x> \}I \\ & = & \{<x \vert C_{g^{-1}} C_g y> + <x \vert A_{g^{-1}} A_g y> \} I \\ & = & < x \vert y> I.
\end{eqnarray*}
\end{proof} 
\medbreak 
We remark further from Theorem \ref{CA} with $g \in {\rm Sp}(V)$ replaced by its inverse that if $v \in V$ then 
	\[c(v) = c_g (C_{g^{-1}} v) + a_g (A_{g^{-1}} v)
\]
	\[a(v) = a_g (C_{g^{-1}} v) + c_g (A_{g^{-1}} v).
\]
\medbreak
\medbreak
Now, the generalized Fock representation of $V$ is set up as follows. For $v \in V$ we define $\pi (v)$ as a complex-linear endomorphism of either the symmetric algebra $SV$ or its full antidual $SV'$ by the rule 
\begin{equation} \label{Fock}
\pi (v) = \frac{1}{\sqrt{2}} \{ c(v) + a(v) \}
\end{equation}
\medbreak
\noindent
whence if $\Phi \in SV'$ and $\psi \in SV$ then 
	\[[\pi (v) \Phi ] (\psi) = \Phi ( \pi (v) \psi).
\]
\medbreak \noindent
Note that if $v \in V$ then as $c(Jv) = i c(v)$ and $a(Jv) = -i a(v)$ so 
	\[c(v) = \frac{1}{\sqrt{2}} \{ \pi (v) - i \pi(Jv) \} 
\]
	\[a(v) = \frac{1}{\sqrt{2}} \{ \pi (v) + i \pi(Jv) \}.
\]
\medbreak 
\medbreak
The generalized Fock representation $\pi$ of $V$ is projective: it satisfies the Heisenberg form of the canonical commutation relations on $SV$ and $SV'$ (without qualification) as follows. 

\medbreak 
\begin{theorem} \label{Heis} 
If $x,y \in V$ then 
	\[[ \pi(x) , \pi(y) ] = i \Omega (x,y) I.
\]
\end{theorem} 
\begin{proof} 
That the displayed equations hold on both $SV$ and $SV'$ follows at once from the canonical commutation relations in Theorem \ref{CCR} : 
\begin{eqnarray*}
	[\pi(x),\pi(y)] & = & \frac{1}{2}[a(x),c(y)] + \frac{1}{2}[c(x),a(y)]\\ & = & \frac{1}{2}\{ <x \vert y> - <y \vert x> \} I \\ & = & i \Omega (x,y) I.
\end{eqnarray*}
\end{proof}
\medbreak 
The generalized Fock representation is also weakly irreducible. 

\medbreak 
\begin{theorem} 
If the linear map $T: V \ra V'$ commutes with $\pi$ in the sense 
	\[v \in V \Rightarrow T \pi (v) = \pi(v) T 
\]
\medbreak 
\noindent 
then $T$ is a scalar (multiple of the canonical inclusion). 
\end{theorem}
\begin{proof}
Here, $\pi(v) \in {\rm End}\: SV$ on the left and $\pi(v) \in {\rm End}\: SV'$ on the right. Taking complex-linear and antilinear parts in the hypothesized condition, if $v \in V$ then $T c(v) = c(v) T$ and $T a(v) = a(v) T$. Now 
	\[v \in V \Rightarrow a(v) T 1 = T a(v) 1 = 0 
\]
\medbreak 
\noindent 
whence Theorem \ref{vac} yields $\lambda \in {\mathbb C}$ such that $T 1 = \lambda 1$. Finally, if $v_1 , \dots ,v_n \in V$ then 
\begin{eqnarray*}
	T(v_1 \cdots v_n) & = & T c(v_1) \cdots c(v_n) 1 \\ & = & c(v_1) \cdots c(v_n) 1 \\ & = & \lambda v_1 \cdots v_n 
\end{eqnarray*}
\medbreak 
\noindent 
and linearity concludes the argument. 
\end{proof} 
\medbreak 
Now, let $g \in {\rm Sp}(V)$. The transformed representation $\pi \circ g$ of $V$ on $SV'$ given by 
\begin{equation} \label{trans}
	v \in V \Rightarrow \pi \circ g (v) = \pi (g v) = \frac{1}{\sqrt{2}} \{ c_g (v) + a_g (v) \} \end{equation} 

\medbreak
\noindent
also satisfies the Heisenberg form of the canonical commutation relations: this may be seen by applying Theorem \ref{CCR'} (rather than Theorem \ref{CCR}) in the proof of Theorem \ref{Heis}. Accordingly, it is reasonable to ask whether the representations $\pi \circ g$ and $\pi$ are equivalent in any sense. 
\medbreak 
By a generalized Fock implementer for $g \in {\rm Sp}(V)$ we shall mean a (nonzero) linear map $U : SV \ra SV'$ that intertwines $\pi$ and $\pi \circ g$ in the sense 
\begin{equation} \label{genimp} 
	v \in V \Rightarrow U \pi (v) = \pi (g v) U 
\end{equation}
\medbreak 
\noindent 
where $\pi(v) \in {\rm End}\: SV$ and $\pi(g v) \in {\rm End}\: SV'$. 
\medbreak 
\begin{theorem} \label{Uca}
The linear map $U : SV \ra SV'$ is a generalized Fock implementer for $g \in {\rm Sp}(V)$ precisely when 
	\[v \in V \Rightarrow \begin{cases} U c(v) = c_g(v) U \\ U a(v) = a_g(v) U. \end{cases}
\]
\end{theorem}
\begin{proof}
In the one direction, taking complex-linear and antilinear parts in the equation \eqref{genimp} defining $U$ as a generalized Fock implementer yields the displayed equations; in the other direction, adding the displayed equations reveals $U$ as a generalized Fock implementer in view of \eqref{Fock} and \eqref{trans}.
\end{proof}
\medbreak

It follows easily by the observation after Theorem \ref{CCR'} that $U : SV \ra SV'$ is a generalized Fock implementer for $g \in {\rm Sp}(V)$ exactly when 
	\[v \in V \Rightarrow \begin{cases} U c_{g^{-1}}(v) = c(v) U \\ U a_{g^{-1}}(v) = a(v) U. \end{cases}
\]
\medbreak
\medbreak
By a generalized Fock vacuum for $g \in {\rm Sp}(V)$ we shall mean a (nonzero) vector $\Phi \in SV'$ such that 
\begin{equation} \label{genvac}
v \in V \Rightarrow \{ \pi(g v) + i \pi(g J v) \} \Phi = 0 
\end{equation}
\medbreak
\noindent
or equivalently 
	\[v \in V \Rightarrow a_g (v)\Phi = 0.
\]
\medbreak
\begin{theorem} \label{UPhi}
If $g \in {\rm Sp}(V)$ then the rule $\Phi = U 1$ sets up a bijective correspondence between the set of all generalized Fock vacua $\Phi \in SV'$ for $g$ and the set of all generalized Fock implementers $U : SV \ra SV'$ for $g$. 
\end{theorem}
\begin{proof}
On the one hand, if $U$ is a generalized Fock implementer and if $v \in V$ then Theorem \ref{Uca} implies that 
	\[a_g (v) U 1 = U a(v) 1 = 0
\]
\medbreak
\noindent
whence $U 1$ is a generalized Fock vacuum. On the other hand, if $\Phi$ is a generalized Fock vacuum then the canonical commutation relations in Theorem \ref{CCR'} enable us to define a generalized Fock implementer $U$ by $U 1 = \Phi$ and the rule that if $v_1, \dots ,v_n \in V$ then 
	\[U(v_1 \cdots v_n) = c_g (v_1) \cdots c_g (v_n) \Phi. 
\]
\medbreak
\noindent
Finally, it is plain that $\Phi \leftrightarrow U$ is a bijective correspondence. 
\end{proof}
\medbreak
Recall that in \eqref{Z} we associated to each $g \in {\rm Sp}(V)$ the symmetric antilinear operator $Z_g = -A_g C_g^{-1}$ with operator norm strictly less than unity; denote the corresponding quadratic by $\zeta_g \in S^2 V'$ so that 
	\[v \in V \Rightarrow a(v) \zeta_g = Z_g v.
\]
\medbreak
\begin{theorem} \label{Gaussvac}
The generalized Fock vacua for $g \in {\rm Sp}(V)$ are precisely the scalar multiples of the Gaussian 
	\[e^{Z_g} = \exp (\zeta_g) \in SV'.
\]
\end{theorem}
\begin{proof}
Let $\Phi = \sum_{d \in {\mathbb N}} \Phi^d$ be a generalized Fock vacuum for $g$. Upon taking homogeneous components, the generalized Fock vacuum condition following \eqref{genvac} on $\Phi$ yields that if $v \in V$ then  $a(C_g v) \Phi^1 = 0$ (the $d = 0$ equation) while 
 \[d > 0 \Rightarrow a(C_g v) \Phi^{d + 1} + c(A_g v) \Phi^{d - 1}.
\]
\medbreak
\noindent
By Theorem \ref{Cinv} it follows that if $v \in V$ then  $a(v) \Phi^1 = 0$ (the $d = 0$ equation) while 
	\[d > 0 \Rightarrow a(v) \Phi^{d + 1} = c(Z_g v) \Phi^{d - 1}.
\]
\medbreak
\noindent
The $d = 0$ equation forces $\Phi^1$ to vanish and the even $d > 0$ equations then force all odd-degree components of $\Phi$ to vanish by Theorem \ref{vac}. The $d = 1$ equation forces $\Phi^2$ to equal $\Phi^0 \zeta_g$ and the odd $d > 0$ equations then force $\Phi = \Phi^0 \exp(\zeta_g)$ by induction. In the opposite direction, each scalar multiple of $\exp(\zeta_g)$ is a generalized Fock vacuum for $g$ either by essentially the same argument or by Theorem \ref{c=a}.
\end{proof}
\medbreak
We are now able to establish the unconditional existence of generalized Fock implementers. 
\medbreak
\begin{theorem} \label{U_g}
The generalized Fock implementers for $g \in {\rm Sp}(V)$ are precisely the scalar multiples of $U_g : SV \ra SV'$ defined by $U_g 1 = e^{Z_g}$ and the rule that if $v_1, \dots ,v_n \in V$ then 
	\[U_g (v_1 \cdots v_n) = c_g (v_1) \cdots c_g (v_n) e^{Z_g}.
\]
\end{theorem}
\begin{proof}
Of course, this is an immediate consequence of Theorem \ref{UPhi} and Theorem \ref{Gaussvac}.
\end{proof}
\medbreak
We remark that if $g \in {\rm Sp}(V)$ then the specific generalized Fock implementer $U_g : SV \ra SV'$ so defined is distinguished by having generalized vacuum expectation value unity in the sense that $[U_g 1](1) = 1$. 
\medbreak
By extension of the usual notion, if $T : SV \ra SV'$ is a linear map then its adjoint is the linear map $T^* : SV \ra SV'$ defined by 
	\[\phi, \psi \in SV \Rightarrow [T^* \phi] (\psi) = \overline{[T \psi] (\phi)}.
\]
\medbreak
\begin{theorem} \label{U*}
If $g \in {\rm Sp}(V)$ then $U_g^* = U_{g^{-1}}$.
\end{theorem}
\begin{proof}
This proceeds with the aid of Theorem \ref{Uca} and the remark thereafter: if $v \in V$ and $\phi, \psi \in SV$ then 
\begin{eqnarray*}
	[U_g^* a(v)\phi](\psi) & = & \overline{[U_g \psi](a(v) \phi)}\\ & = & \overline{[c(v) U_g \psi](\phi)}\\ & = & \overline{[U_g c_{g^{-1}}(v) \psi](\phi)}\\ & = & [U_g^* \phi] (c_{g^{-1}}(v) \psi) \\ & = & [a_{g^{-1}}(v)U_g^* \phi](\psi)
\end{eqnarray*}
\medbreak
\noindent
whence
	\[U_g^* a(v) = a_{g^{-1}} U_g^*
\]
\medbreak
\noindent
while 
	\[U_g^* c(v) = c_{g^{-1}} U_g^*
\]
\medbreak
\noindent
similarly; finally, 
	\[[U_g^* 1] (1) = \overline{[U_g 1] (1)} = 1.
\]
\end{proof}
\medbreak
Now traditionally, the Fock representation and Fock implementers act in symmetric Fock space $S[V]$. The relationships between our generalized notions and the traditional ones are as follows. 
\medbreak
First of all, let $v \in V$. The generalized Fock operator $\pi(v) : SV' \ra SV'$ restricts to define in $S[V]$ an operator also denoted by $\pi(v)$ having natural domain 
	\[\{ \Phi \in S[V] : \pi(v) \Phi \in S[V] \}.
\]
\medbreak
\noindent
An argument along similar lines to that for Theorem \ref{adjoints} shows that this traditional Fock operator $\pi(v)$ with the above domain is self-adjoint: $\pi(v)^* = \pi(v)$. We point out that Theorem \ref{Heis} is not true for these traditional Fock operators without qualification: domain technicalities enter into the (Heisenberg) canonical commutation relations, thus 
	\[x,y \in V \Rightarrow [\pi(x), \pi(y)] \subset i \Omega(x,y) I.
\]
\medbreak
Again let $g \in {\rm Sp}(V)$. In the traditional context, it is natural to seek conditions necessary and sufficient for the existence of a unitary operator $U:S[V] \ra S[V]$ such that 
	\[v \in V \Rightarrow U \pi(v) = \pi(g v) U.
\]
\medbreak
\noindent
As Theorem \ref{U_g} furnishes a linear map $U_g : SV \ra SV'$ such that 
	\[v \in V \Rightarrow U_g \pi(v) = \pi(g v) U_g 
\]
\medbreak
\noindent
it is clear that the problem to solve now is essentially one of normalization. 
\medbreak
\begin{theorem} \label{unitU}
If $g \in {\rm Sp}(V)$ is such that $A_g$ is of Hilbert-Schmidt class then the prescription
 \[ U(g) :\:= \Vert e^{Z_g} \Vert^{-1} U_g
\]
\medbreak
\noindent
determines a unitary operator on $S[V]$. 
\end{theorem}
\begin{proof}
Let $A_g$ be Hilbert-Schmidt. The symmetric antilinear operator $Z_g$ is now Hilbert-Schmidt also; as $\Vert Z_g \Vert < 1$ already, Theorem \ref{Gauss} places $e^{Z_g}$ in $S[V]$ with 
	\[\Vert e^{Z_g} \Vert^4 = {\rm Det} (I - Z_g^2)^{-1}.
\]
\medbreak
\noindent
Normalizing, define $U(g) = \Vert e^{Z_g} \Vert^{-1} U_g$ as announced. The corresponding generalized Fock vacuum $\Phi(g) = U(g) 1 = \Vert e^{Z_g} \Vert^{-1} e^{Z_g} \in S[V]$ is a unit vector in the domain of every creator-annihilator polynomial, on account of the remark after Theorem \ref{smooth}. From the definition of $U_g$ in Theorem \ref{U_g} it now follows that $U(g)$ maps $SV$ to $S[V]$. To see that $U(g): SV \ra S[V]$ is isometric, let $x_1, \dots ,x_r, y_1, \dots ,y_s \in V$:  the canonical commutation relations in Theorem \ref{CCR'} yield 
\begin{eqnarray*}
	<U(g)(x_1 \cdots x_r) \vert U(g)(y_1 \cdots y_s)> & = & <c_g(x_1) \cdots c_g(x_r) \Phi(g) \vert c_g(y_1) \cdots c_g(y_s) \Phi(g)> \\ & = & <\Phi(g) \vert a_g(x_r) \cdots a_g(x_1) c_g(y_1) \cdots c_g(y_s) \Phi(g)> \\ & = & <x_1 \cdots x_r \vert y_1 \cdots y_s>
\end{eqnarray*}
\medbreak
\noindent
by virtue of Theorem \ref{adjoints}. Of course, parallel remarks apply to $U(g^{-1})$ because $Z_{g^{-1}} = C_g^{-1} A_g$ is Hilbert-Schmidt. To see that the isometric extension $\overline{U(g)}: S[V] \ra S[V]$ is unitary, note first that $U(g^{-1}) = U(g)^*$ by Theorem \ref{U*} and the fact that $I - Z_{g^{-1}}^2 = C_g^{-1} (I - Z_g^2) C_g$ from Theorem \ref{CA}. Now, if $\phi, \psi \in SV$ then Theorem \ref{compat} shows that 
	\[<\phi \vert U(g)\psi> = [U(g)\psi](\phi) = \overline{[U(g)^* \phi](\psi)} = <U(g^{-1})\phi \vert \psi> 
\]
\medbreak
\noindent
whence if $\Phi, \Psi \in S[V]$ then Theorem \ref{dense} shows that 
	\[< \Phi \vert \overline{U(g)}\Psi> = <\overline{U(g^{-1})}\Phi \vert \Psi>.
\]
\medbreak
\noindent
Thus the Hilbert space adjoint of $\overline{U(g)}$ is the isometry $\overline{U(g^{-1})}$.

\end{proof}
\medbreak
Conversely, if $U_g$ may be rescaled so as to produce a unitary operator on $S[V]$ then in particular the Gaussian $e^{Z_g} = U_g 1$ lies in $S[V]$ and therefore $A_g = -Z_g C_g$ is of Hilbert-Schmidt class. 
\medbreak
Now by definition, the restricted symplectic group ${\rm Sp}_{\rm res}(V)$ comprises precisely all those $g \in {\rm Sp}(V)$ for which $A_g$ is of Hilbert-Schmidt class. When $g \in {\rm Sp}_{\rm res}(V)$ we shall denote by 
\begin{equation} 
\overline{U}(g) = {\rm Det}^{\frac{1}{4}} (I - Z_g^2) \overline{U_g}
\end{equation}
\medbreak
\noindent
the extension of $U(g) = {\rm Det}^{\frac{1}{4}} (I - Z_g^2) U_g$ to a unitary operator on $S[V]$. By definition, the resulting map 
\begin{equation} 
\overline{U} : {\rm Sp}_{\rm res}(V) \ra {\rm Aut}S[V]
\end{equation}
\medbreak
\noindent
is the metaplectic representation. This is indeed a projective representation, whose cocycle may be derived explicitly as follows. 
\medbreak
\begin{theorem}
If $g, h \in {\rm Sp}_{\rm res}(V)$ then 
	\[\overline{U_g}\ \overline{U_h} = \delta (g,h)\overline{U_{gh}}
\]
\medbreak
\noindent
where 
	\[\delta (g,h) = {\rm Det}^{\frac{1}{2}} (I - Z_h Z_{g^{-1}})^{-1}.
\]
\end{theorem} 
\begin{proof}
Introduce a linear map $\widetilde{U}_{gh} :SV \ra SV'$ by the rule 
	\[\phi \in SV \Rightarrow \widetilde{U}_{gh}(\phi) = \overline{U_g}\ \overline{U_h}(\phi) = \overline{U_g}(U_h \phi).
\]
\medbreak
\noindent
If $v \in V$ then it follows by Theorem \ref{compat} and Theorem \ref{adjoints} with the proof of Theorem \ref{unitU} that 
\begin{eqnarray*}
	[\widetilde{U}_{gh} c(v) \phi](\psi) & = & <\psi \vert \overline{U}_g U_h c(v) \phi> = <U_{g^{-1}} \psi \vert c_h (v) U_h \phi > \\ & = & <a_h (v) U_{g^{-1}} \psi \vert U_h \phi> = <U_{g^{-1}} a_{gh} (v) \psi \vert U_h \phi> \\ & = & <a_{gh} (v) \psi \vert \overline{U}_g U_h \phi> = [\widetilde{U}_{gh} \phi](a_{gh}(v) \psi) \\ & = & [c_{gh} (v) \widetilde{U}_{gh} \phi](\psi).
\end{eqnarray*}
\medbreak
\noindent
Accordingly, if $v \in V$ then 
	\[\widetilde{U}_{gh} c(v) = c_{gh}(v) \widetilde{U}_{gh}
\]
\medbreak
\noindent
and similarly 
	\[\widetilde{U}_{gh} a(v) = a_{gh}(v) \widetilde{U}_{gh}.
\]
\medbreak
\noindent
Thus Theorem \ref{Uca} and Theorem \ref{U_g} imply that $\widetilde{U}_{gh}$ and $U_{gh}$ are proportional, so $\overline{U_g}\ \overline{U_h}$ and $\overline{U_{gh}}$ are proportional. All that remains is to compare normalizations: on the one hand, $[\overline{U_{gh}} 1](1) = 1$ by definition; on the other hand, Theorem \ref{compat} and Theorem \ref{pair} with the proof of Theorem \ref{unitU} yield 
\begin{eqnarray*}
	[\overline{U_g}\ \overline{U_h} 1](1) & = & <1 \vert \overline{U_g}(U_h 1)> = <U_{g^{-1}} 1 \vert U_h 1> \\ & = & <e^{Z_{g^{-1}}} \vert e^{Z_h}> = {\rm Det}^{\frac{1}{2}} (I - Z_h Z_{g^{-1}})^{-1}.
\end{eqnarray*}
\end{proof}

\section{Remarks}
\medbreak
In this final section, we make a number of remarks concerning the approach adopted in these notes. 
\medbreak
Firstly, the approach via the antidual is decidedly elegant and offers a natural environment in which to develop the theory. It facilitates clean proofs: indeed, we have taken this opportunity to present simple proofs for several theorems difficult to locate in the literature. Thus, the handling of creators and annihilators is improved: for example, the proofs that if $v \in V$ then $c(v)^* = a(v)$ and $a(v)^* = c(v)$ are particularly straightforward; field operators and the number operator are similarly transparent. Also, exponentials are manipulated with ease: among other things, we mention the effect of creators and annihilators on Gaussians and the fact that the exponentials of nonzero cubics do not lie in symmetric Fock space. Of course, the antidual is especially appropriate for the discussion of generalized Fock implementation. 
\medbreak
As another example, let us outline a proof of the fact that if $Z \in \Sigma^2 [V]$ and $\Vert Z \Vert < 1$ then the Gaussian $e^Z \in S[V]$ is cyclic for creators alone. Observe that $I - Z^2$ is an invertible positive operator, so we may define $C: \: = \sqrt{(I - Z^2)^{-1}}$; the operator $g: \: = (I - Z)C$ then lies in ${\rm Sp}(V)$ and indeed in ${\rm Sp}_{\rm res}(V)$ since $Z_g = Z$ is Hilbert-Schmidt. Now, the unitary operator $\overline{U(g)}$ on $S[V]$ defined in Theorem \ref{unitU} has the property that if $v_1, \dots ,v_n \in V$ then 
	\[\overline{U(g)}(v_1 \cdots v_n) = \Vert e^Z \Vert^{-1} c_g (v_1) \cdots c_g (v_n) e^Z
\]
\medbreak
\noindent
whence Theorem \ref{c=a} implies that 
	\[\overline{U(g)}(v_1 \cdots v_n) \in \{ \phi \:e^Z : \phi \in SV \}.
\]
\medbreak
\noindent
As the (possibly empty) products of vectors from $V$ span $SV$ and as $\overline{U(g)}$ is unitary, so $\{\phi \:e^Z : \phi \in SV \}$ is dense in $S[V]$. Otherwise said, $e^Z$ is cyclic for creators alone. 
\medbreak
Next, we ought at least to mention the direct construction of the bosonic Fock representation in Weyl form. Coherent states are especially well-suited for this purpose, so let us introduce a complex vector space $EV$ with basis $\{ \varepsilon^z : z \in V \}$ and inner product given by the rule that if $x,y \in V$ then $ < \varepsilon^x \vert \varepsilon^y > = e^{<x \vert y>}$. Notice that Theorem \ref{enorm} and Theorem \ref{total} permit us to identify $EV$ with the span of the coherent vectors $\{ e^z : z \in V \}$. Along with $EV$ itself we naturally consider its full antidual $EV'$ whose subspace $E[V]$ of bounded antilinear functionals on $EV$ is identified with $S[V]$. Certain other subspaces of $EV'$ are also important: for example, that comprising all $\Phi \in EV'$ for which the function $V \ra {\mathbb C} : z \mapsto \Phi(\varepsilon^z)$ is antiholomorphic in one of several senses, such as the usual sense on finite-dimensional subspaces. 
\medbreak
To each $v \in V$ we associate the linear automorphism $W(v)$ of $EV$ defined by the rule  
	\[z \in V \Rightarrow W(v)\varepsilon^z = (\Vert \varepsilon^v \Vert e^{<v \vert z>})^{-1} \varepsilon^{v + z} 
\]
\medbreak 
\noindent 
and extend it to $EV'$ by antiduality according to the prescription 
	\[\Phi \in EV', \psi \in EV \Rightarrow [W(v)\Phi](\psi) = \Phi(W(-v)\psi). 
\]
\medbreak
\noindent
Direct computation reveals that $W(v)$ is unitary on $EV$ and indeed on $E[V]$. The resulting map $W : V \ra {\rm Aut} E[V]$ is a regular projective representation: it is regular, for if $x, y, v \in V$ then the inner product 
	\[<\varepsilon^x \vert W(t v) \varepsilon^y> = \exp {\{ <x \vert y> + (<x \vert v> - <v \vert y>)t - \frac{1}{2}\Vert v \Vert^2 t^2 \}}
\]
\medbreak
\noindent
depends continuously on $t \in {\mathbb R}$; it is projective, its cocycle being readily verified to have the Weyl form 
	\[x, y \in V \Rightarrow W(x) W(y) = \exp \{- i \Omega(x,y)\} W(x + y).
\]
\medbreak 

In this formalism, a generalized Fock implementer for $g \in {\rm Sp}(V)$ is a (nonzero) linear map $U: EV \ra EV'$ that intertwines $W$ on $EV$ with $W \circ g$ on $EV'$ in the sense 
	\[v \in V \Rightarrow U W(v) = W(g v) U.
\]
\medbreak
\noindent
The intertwiner $U$ may be required to satisfy further restrictions, such as that $<\varepsilon^x \vert U \varepsilon^y>$ be (antiholomorphic, holomorphic) in $(x, y) \in V \times V$. With this definition, a specific generalized Fock implementer $U_g : EV \ra EV'$ is given explicitly by the rule that if $x, y \in V$ then 
	\[[U_g \varepsilon^y](\varepsilon^x) = \exp{\{\frac{1}{2}<x \vert C_{g^{-1}}^{-1} (y - A_{g^{-1}} x)> + \frac{1}{2}<C_g ^{-1} (x - A_g y) \vert y> \}}
\]
\medbreak
\noindent
The proof of this fact is entirely routine: as the action of $W$ passes from $EV$ to $EV'$ by antiduality, it is enough to argue algebraically that if $x, y, v \in V$ then 
	\[[U_g W(v) \varepsilon^y](\varepsilon^x) = [U_g \varepsilon^y](W(-g v) \varepsilon^x). 
\]
\medbreak 
\noindent
Of course, if $g \in {\rm Sp}_{\rm res}(V)$ then ${\rm Det}^{\frac{1}{4}} (I - Z_g^2) U_g$ determines a unitary intertwining operator on $E[V]$. We remark that \cite{PSSZ} presents a more detailed analysis, incorporating $(- , +)$  holomorphicity restrictions in terms of the complex-wave representation. 
\medbreak
Lastly, the elegance of the approach adopted here suggests that it should be adopted elsewhere. As a matter of fact, in \cite{RIMS} we have already discussed an analogous treatment for the fermionic Fock representation of $V$: we placed fermionic Fock space $\bigwedge[V]$ between the exterior algebra $\bigwedge V$ and its full antidual $\bigwedge V'$ while simultaneously developing the Berezin calculus in arbitrary dimensions. In the fermionic context, it transpires that an orthogonal transformation $g \in O(V)$ admits a generalized Fock implementer precisely when the complex-linear part $C_g$ has finite-dimensional kernel; again, if the antilinear part $A_g$ is Hilbert-Schmidt then a suitably normalized implementer determines a unitary intertwining operator on $\bigwedge [V]$. Of course, it is natural to attempt a similar treatment for the Fock representation of an indefinite inner product space : when this is a Krein space the Hilbert space machinery may be employed, but even then it is not of primary importance; thus an approach by way of the antidual shows promise. Such matters will be addressed in a future publication.

\medbreak


\begin{thebibliography} {99}

\bibitem{Araki} 
H. ARAKI, On Quasifree States of the Canonical Commutation Relations (II). \emph{Publ. RIMS, Kyoto Univ.} \textbf{7} (1971/72) 121-152. 

\bibitem{BSZ} 
J. C. BAEZ, I. E. SEGAL and Z. ZHOU, \emph{Introduction to Algebraic and Constructive Quantum Field Theory}. Princeton University Press (1992). 

\bibitem{Ber} 
F. A. BEREZIN, \emph{The Method of Second Quantization}. Academic Press (1966). 

\bibitem{BraRob} 
O. BRATTELI and D. W. ROBINSON, \emph{Operator Algebras and Quantum Statistical Mechanics II}. Springer-Verlag (1981). 

\bibitem{HY}
Z. Y. HUANG and J. A. YAN, \emph{Introduction to Infinite Dimensional Stochastic Analysis}. Kluwer Academic (2000). 

\bibitem{Obata} 
N. OBATA, \emph{White Noise Calculus and Fock Space}. Springer-Verlag (1994). 


\bibitem{PSSZ} 
S. M. PANEITZ, J. PEDERSEN, I. E. SEGAL and Z. ZHOU, Singular Operators on Boson Fields as Forms on Spaces of Entire Functions on Hilbert Space. \emph{J. Functional Analysis} \textbf{100} (1991) 36-58. 

\bibitem{RIMS} 
P.L. ROBINSON, The Berezin Calculus. \emph{Publ. RIMS, Kyoto Univ.} \textbf{35} (1999) 123-194. 

\bibitem{Rui} 
S. N. M. RUIJSENAARS, On Bogoliubov transformations. II. The general case. \emph{Annals of Physics} \textbf{116} (1978) 105-134. 

\bibitem{GBS} 
G. B. SEGAL, Unitary Representations of some Infinite Dimensional Groups. \emph{Comm. Math. Phys.} \textbf{80} (1981) 301-342. 

\bibitem{Shale} 
D. SHALE, Linear symmetries of free boson fields. \emph{Trans. Amer. Math. Soc.} \textbf{103} (1962) 149-167. 

\bibitem{Vergne} 
M. VERGNE, Groupe symplectique et seconde quantification. \emph{C. R. Acad. Sci. S\'er. A} \textbf{285} (1977) 191-194. 



\end{thebibliography}
\end{document}